\documentclass[]{article}  
\usepackage{CJK}
\usepackage{amssymb}
\usepackage{bbm}
\usepackage{amsfonts}
\usepackage{mathrsfs}
\usepackage{graphicx}
\usepackage{amsmath}
\usepackage{amsthm}
\usepackage{xypic}
\usepackage{graphicx}
\usepackage{times}
\usepackage{geometry}
\usepackage{color}
\usepackage{indentfirst}
\usepackage{cases}
\usepackage{hyperref}
\usepackage{graphicx}
\usepackage{float}
\usepackage{amsthm}
\usepackage{amssymb}
\usepackage{amsmath}
\usepackage{mathrsfs}
\usepackage{indentfirst}
\usepackage{geometry}
\usepackage{authblk}
\usepackage{hyperref}
\usepackage{enumerate}
\usepackage[numbers,sort&compress]{natbib}
\usepackage{amsmath, amssymb}
\usepackage{indentfirst}
\usepackage{geometry}
\usepackage{authblk}  
\usepackage{hyperref}
\usepackage{enumerate}
\usepackage{tikz}
\usetikzlibrary{arrows.meta}
\usepackage{caption}
\usepackage{subcaption}

\newtheorem{Theorem}{Theorem}[section]
\theoremstyle{definition}
\newtheorem{Corollary}[Theorem]{Corollary}

\newtheorem{Lemma}[Theorem]{Lemma}

\newtheorem{Prop}{Proposition}

\theoremstyle{definition}

\title{Some Notes on 2-PDR of Valency 3}
\author{Songnian Xu \thanks{Corresponding author. E-mail address: xsn131819@163.com},\ \ \ Dein Wong \thanks{Corresponding author. E-mail address:wongdein@163.com.  Supported by the National Natural Science Foundation of China(No.12371025)} ,\ \ Chi Zhang \thanks{ Supported by NSFC of China(No.12001526)}, \ \ Wenhao Zhen
\\ {\small  \it School of Mathematics, China University of Mining and Technology, Xuzhou,  China.}

}
\date{}

\begin{document}
\baselineskip 17pt

\title{$m$-partite oriented semiregular representation of valency 3 for finite groups}

 \author{
Songnian Xu\textsuperscript{a,}\thanks{Corresponding author. E-mail address: xsn131819@163.com},
\ Dein Wong.\textsuperscript{a,}\thanks{Corresponding author. E-mail address: wongdein@163.com},
\ Wenhao Zhen\textsuperscript{a}
\\
\textsuperscript{a}School of Mathematics, China University of Mining and Technology, Xuzhou, China.

}
\date{}
\maketitle

\begin{abstract}
Let $G$ be a finite group and $m \geq 2$ a positive integer. We say that $G$ admits an \emph{oriented $m$-semiregular representation} (abbreviated as OmSR) if there exists a  $m$-Cayley digraph $\Gamma$ over $G$ such that $\Gamma$ is oriented and $\mathrm{Aut}(\Gamma) \cong G$.
In \cite{xu1}, we classified finite groups generated by at most two elements that admit an OmSR of valency 3 for $m \geq 2$. In this article, we consider $m$-partite digraphs.We say a finite group $G$ admits an \emph{$m$-partite oriented semiregular representation} ($m$-partite digraphical representation), abbreviated as \emph{$m$-POSR} (\emph{$m$-PDR}), if there exists an \emph{oriented} $m$-partite Cayley digraph (\emph{$m$-partite Cayley digraph}) $\Gamma$ with $\mathrm{Aut}(\Gamma) \cong G$.
In this paper, we classify finite groups generated by at most two elements that admit $m$-POSR. Since if $G$ admits an $m$-POSR, then $G$ must also admit an $m$-PDR (while the converse does not hold), as a natural consequence, we also provide a complete classification for groups $G=\langle x,y\rangle$ that admit $m$-PDR of valency 3. This complements the results in \cite{xu2}.
\end{abstract}

\let\thefootnoteorig\thefootnote
\renewcommand{\thefootnote}{\empty}

\section{Introduction}
A \textit{directed graph} (or \textit{digraph}) $\Gamma$ is an ordered pair $(V(\Gamma), A(\Gamma))$, where $V(\Gamma)$ is a non-empty vertex set and $A(\Gamma) \subseteq V(\Gamma) \times V(\Gamma)$. Elements of $V(\Gamma)$ and $A(\Gamma)$ are called \textit{vertices} and \textit{arcs}, respectively. For an arc $(u, v) \in A(\Gamma)$, $v$ is an \textit{out-neighbor} of $u$, and $u$ is an \textit{in-neighbor} of $v$.
The $out$-$valency$ (resp. $in$-$valency$) of a vertex $v \in V(\Gamma)$ counts its out-neighbors (resp. in-neighbors). A digraph $\Gamma$ is \textit{$k$-regular} if some positive integer $k$ ensures every vertex has equal out-degree and in-degree $k$.
In this paper, we require that all digraphs are regular.

For $X \subseteq V(\Gamma)$, the \textit{induced subdigraph} on $X$ is $\Gamma[X] := (X, A(\Gamma) \cap (X \times X))$, abbreviated as $[X]$. A digraph is a \textit{graph} if $A(\Gamma)$ is symmetric (i.e., $(u, v) \in A(\Gamma)$ implies $(v, u) \in A(\Gamma)$), and \textit{oriented} if for distinct $u, v$, at most one of $(u, v)$ or $(v, u)$ lies in $A(\Gamma)$.

An \textit{automorphism} of $\Gamma$ is a permutation $\sigma$ of $V(\Gamma)$ preserving arcs: $(x^\sigma, y^\sigma) \in A(\Gamma) \iff (x, y) \in A(\Gamma)$. The set of all such automorphisms forms the \textit{full automorphism group} $\operatorname{Aut}(\Gamma)$.

Let \( G \) be a permutation group acting on a set \( X \). For any \( x \in X \), let \( G_x \) stand for the stabilizer of \( x \) within \( G \). We say that \( G \) acts \textit{semiregularly} on \( X \) if \( G_x = \{e\} \) for every \( x \in X \); it acts \textit{regularly} if it is both semiregular and transitive.
It should be noted that the regularity of a group \( G \) acting on a set \( X \) and the regularity of a graph \( \Gamma \) are two distinct concepts.

Let \( m \) be a positive integer. An \textit{\( m \)-Cayley digraph} of a finite group \( G \) is defined as a digraph \( \Gamma \) that has a semiregular automorphism group isomorphic to \( G \), with this group having exactly \( m \) orbits on the vertex set of \( \Gamma \). In particular, when \( m = 1 \), \( \Gamma \) is precisely a \textit{Cayley digraph}.

A classic question put forward by K$\ddot{o}$nig \cite{kon} in 1936 is: Can a given group be represented as the automorphism group of some graph? K$\ddot{o}$nig's work \cite{kon} is regarded as the starting point of group representation theory. This problem was solved by Frucht \cite{fru} in 1939, and he further refined the result in \cite{fru} by showing that every finite group can be the automorphism group of a cubic graph.

A group \( G \) is said to have a \textit{digraphical regular representation} (DRR) or a \textit{graphical regular representation} (GRR) if there exists a digraph or graph \( \Gamma \), respectively, such that \( \mathrm{Aut}(\Gamma) \cong G \) and acts regularly on the vertex set of \( \Gamma \). One of the most fundamental questions in this area is the "GRR and DRR problem": \textit{Which groups have a GRR or a DRR?}

Babai (\cite{bab1}, Theorem 2.1) showed that all groups admit a DRR, with the only exceptions being \( Q_8, \mathbb{Z}_2^2, \mathbb{Z}_2^3, \mathbb{Z}_2^4, \) and \( \mathbb{Z}_3^2 \). Morris and Spiga \cite{mor,mor1,spi1} provided a classification of finite groups that admit an ORR. It is obvious that if a group has a GRR, then it must also have a DRR, but the reverse is not true. After a long sequence of partial results from various researchers\cite{he,im,im1,im2,no,no1,wa}, Godsil completed the classification of groups with a GRR in \cite{god1}.

A group \( G \) is said to admit a \textit{graphical \( m \)-semiregular representation} (GmSR), \textit{digraphical \( m \)-semiregular representation} (DmSR), or \textit{oriented \( m \)-semiregular representation} (OmSR) if there exists an \( m \)-Cayley graph, \( m \)-Cayley digraph, or oriented \( m \)-Cayley  digraph \( \Gamma \) of $G$, respectively, such that \( \mathrm{Aut}(\Gamma) \cong G \). Finite groups admitting a GmSR or DmSR have been classified by Du and coauthors in \cite{du1,du1'}.

The classification problems for finite groups admitting GmSR, OmSR, DmSR, GRR, ORR, and DRR of valency $k$ remain largely open for any fixed integer $k \geq 2$. Verret and Xia \cite{x2} classified finite simple groups with an ORR (i.e., O1SR) of valency 2, proving that every simple group of order at least 5 has such an ORR.
Du \cite{du2} classified 2-generated groups with an OmSR of valency 2, while the authors \cite{xu1} addressed 2-generated groups with an OmSR of valency 3.

In this paper, we focus on \( m \)-partite digraphs. An \textit{\( m \)-partite digraph} is a digraph whose vertex set can be divided into \( m \) subsets (or parts) such that there are no arcs within any single part (equivalently, the subgraph induced by each part is an empty digraph).
We say that a finite group \( G \) admits an \textit{\( m \)-partite oriented semiregular representation} ($m$-POSR) if there exists an oriented \( m \)-partite  Cayley digraph \( \Gamma \) over \( G \) with \( \mathrm{Aut}(\Gamma) \cong G \). For \( m \)-partite Cayley digraphs, it is easy to see that a 1-partite Cayley digraph must be an empty graph. Therefore, throughout this paper, we only consider the case where \( m \geq 2 \).

\begin{Theorem}
Let $G = \langle x \rangle$ be a finite cyclic group and $m \geq 2$ be a positive integer. Then $G$ admits an $m$-POSR of valency 3, except when either:
\begin{enumerate}

\item[(i)]  $m=2$ and $o(x)\leq 6 $,
\item[(ii)]  $m=3$ and $o(x)\leq 3$,
\item[(iii)] $m=4$ and $o(x)\leq2$.
\item[(iv)] $5\leq m\leq 8$ and $o(x)=1$
\end{enumerate}
\end{Theorem}

Let $\Phi = \{ G \mid G = \langle x, y \rangle, o(x) = 4 \text{ and for any } \langle a, b \rangle = G, o(a), o(b) \leq 4 \}$.

\begin{Theorem}
    Let $G = \langle x, y \rangle \neq \langle x \rangle$ and $G \in \Phi$. Then $G$ admits a 2-POSR of valency 3 except for $G\in \{Q_8, C_4\rtimes C_4, (C_2\times C_2)\rtimes C_4, (C_4\times C_2)\rtimes C_4$\}.
\end{Theorem}

\begin{Theorem}
Let $G = \langle x, y \rangle \neq \langle x \rangle$ and $G \notin \Phi$. Then $G$ admits a 2-POSR of valency 3 except for the following cases:
\[
G \in \{ \mathbb{Z}_2^2, D_6 \},
\]
where $D_n$ denotes the dihedral group of order $n$.
\end{Theorem}

\begin{Theorem}
Let $G = \langle x, y \rangle\neq \langle x\rangle$ be a finite  group. Then for every integer $m \geq 3$, $G$ admits an $m$-POSR of valency 3.
\end{Theorem}

Theorems 1.1, 1.2, 1.3 and 1.4 have provided a complete classification for groups $G=\langle x,y\rangle$ that admit an $m$-POSR of valency 3 with $m\geq 2$. Since any group admitting an $m$-POSR must also admit an $m$-PDR, we obtain the following corollary after a straightforward verification. This completes the results in \cite{xu2}.

\begin{Theorem}\cite{xu2}
Let $G = \langle x \rangle \neq 1$ be a cyclic group. Then for every integer $m\geq2$, $G$ admits an $m$-PDR of valency 3, except for when:
\begin{enumerate}
    \item[(i)] $m = 2$ and $o(x) \leq 4$,
    \item[(ii)] $m = 3$ and $o(x)=2$ .
\end{enumerate}

\end{Theorem}

\begin{Corollary}
    Let $G=\langle x,y\rangle $ be a finite group and $m\geq 2$ an integer. Then $G$ admits an $m$-PDR of valency 3 except for the following cases:
    
\noindent (1) $G=<x>$ and: 
\begin{enumerate}   
\item[(i)]  $m=2$ and $o(x)\leq 4 $,
\item[(ii)]  $m=3$ and $o(x)\leq 2$,
\item[(iii)] $3\leq m\leq 6$ and $o(x)=1$.
\end{enumerate}   
\noindent (2) $G=\mathbb{Z}^2_2$ and $m=2$

\end{Corollary}

\section{Preliminaries and notations}

For a positive integer \( m \) and finite group \( G \), we denote \( (g,i) \in G \times \{0,\ldots,m-1\} \) as \( g_i \) for brevity, with \( \{0,\ldots,m-1\} \) identified with \( \mathbb{Z}_m \) (integers modulo \( m \)). Recall that an \( m \)-Cayley digraph of \( G \) is a digraph with a semiregular automorphism group (isomorphic to \( G \)) having exactly \( m \) vertex orbits. We now give a more precise definition of this structure.

For each \( i \in \mathbb{Z}_m \) and \( G_i= \{g_i \mid g \in G\} \). Similar to classical Cayley digraphs, an \( m \)-Cayley digraph may be represented as
\[
\Gamma = \mathrm{Cay}(G, T_{i,j} : i, j \in \mathbb{Z}_m)
\]
with:
\begin{itemize}
\item Vertex set: \( G \times \mathbb{Z}_m = \bigcup_{i \in \mathbb{Z}_m} G_i \),
\item Arc set: \( \bigcup_{i,j \in \mathbb{Z}_m} \{(g_i, (tg)_j) \mid t \in T_{i,j}\} \),
\end{itemize}
The digraph $\Gamma$ is called:
\begin{itemize}
    \item An \emph{oriented $m$-Cayley digraph} if $T_{i,j}\cap T^{-1}_{j,i} = \emptyset$ for all $i,j \in \mathbb{Z}_m$,
    \item An \emph{$m$-partite Cayley digraph} if $T_{i,i} = \emptyset$ for every $i \in \mathbb{Z}_m$.
\end{itemize}
For any \( g \in G \), the right multiplication map \( R(g) \)—which acts by \( R(g): x_i \mapsto (xg)_i \) for all \( x_i \in G_i \) and \( i \in \mathbb{Z}_m \)—is an automorphism of \( \Gamma \). The set \( R(G) = \{R(g) \mid g \in G\} \) is isomorphic to \( G \) and forms a semiregular automorphism group of \( \Gamma \) with orbits given by the \( G_i \).

For a digraph \( \Gamma \) and any vertex \( x \in V(\Gamma) \), we define:
\begin{itemize}
  \item \( \Gamma^{+0}(x) = \{x\} \),
  \item \( \Gamma^{+1}(x) = \Gamma^+(x) \) (the set of out-neighbors of \( x \) in \( \Gamma \)) and \( \Gamma^+[x] = \Gamma^+(x) \cup \{x\} \),
  \item \( \Gamma^{+2}(x) \) as the union of out-neighbors of all vertices in \( \Gamma^+(x) \) and ,
  \item \( \Gamma^{+k}(x) = \bigcup_{y \in \Gamma^{+(k-1)}(x)} \Gamma^+(y) \) for any integer \( k \geq 1 \).
\end{itemize}

\begin{Prop}\cite[3.1]{du3}
Let $m$ be a positive integer at least $2$ and let $G$ be a finite group.
For any $i, j\in \mathbb{Z}_m$, let $T_{i,j}\subseteq G$ and let $\Gamma=Cay(G,T_{i,j}: i,j\in \mathbb{Z}_m)$ be a connected $m$-cayley digraph over $G$.
For $A=Aut(\Gamma)$, if $A$ fixes $G_i$ setwise for all $i\in \mathbb{Z}_m$ and there exist $u_0\in G_0,u_1\in G_1,\ldots,u_{m-1}\in G_{m-1}$ such that $A_{u_i}$ fixes $\Gamma^{+}(u_i)$ pointwise for all $i\in \mathbb{Z}_m$, then $A=R(G)$.
\end{Prop}

\section{Proof of Theorems}
It is worth noting that, according to the definitions of OmSR and \( m \)-POSR, if a group \( G \) admits an \( m \)-POSR, then \( G \) must necessarily admit an OmSR, though the converse is not true. However, when \( G \cong \mathbb{Z}_2 \) or \( \mathbb{Z}_2^2 \), saying that \( G \) admits an OmSR is equivalent to saying that \( G \) admits an \( m \)-POSR. 
This is because, for an oriented \( m \)-Cayley digraph \( \Gamma = (G, T_{i,j} : i \in \mathbb{Z}_m) \), when \( G \cong \mathbb{Z}_2 \) or \( \mathbb{Z}_2^2 \), we have \( T_{i,i} = \emptyset \) since \( T_{i,i} \cap T_{i,i}^{-1} = \emptyset \), which means \( \Gamma \) is an oriented \( m \)-partite Cayley digraph in this case. 

Moreover, in \cite{xu1}, for \( G \cong \mathbb{Z}_2 \) or \( \mathbb{Z}_2^2 \), we have completed the classification of \( G \) admitting an OmSR of valency 3, with the results given in Lemma 3.1 and Lemma 3.2 below. Thus, throughout the rest of this paper, we always require that $o(x) \geq 3$ or $o(y) \geq 3$ when $G = \langle x, y \rangle$, where \( o(x) \) denotes the order of \( x \).

\begin{Lemma}\cite[3.1]{xu1}
Let $G=\langle x\rangle\cong \mathbb{Z}_2$ be a finite group and $m\geq 2$ be a positive integer.
Then $G$ admits an OmSR of valency 3 if and only if $m\geq 5$.
\end{Lemma}

\begin{Lemma}\cite[3.3]{xu1}
Let $G=\langle x,y\rangle\cong \mathbb{Z}_2^2$ and $m\geq 2$ an integer. Then $G$ admits an OmSR of valency 3, except when $m=2$.
\end{Lemma}

When $G = 1$, $G$ admits an $m$-POSR of valency 3 if and only if there exists a 3-regular antisymmetric oriented digraph of order $m$, where an antisymmetric (di)graph is defined as a (di)graph whose automorphism group is trivial.
 \begin{Lemma}
Let $m$ be a positive integer. There exists a 3-regular antisymmetric  oriented digraph of order $m$ if and only if $m \geq 9$.
\end{Lemma}
\begin{proof}
    Using Mathematica, we find that there exists no 3-regular antisymmetric oriented digraph of order $m$ for $m \leq 8$.

\begin{figure}[h]
\centering

\begin{subfigure}[b]{0.45\textwidth}
\centering
\begin{tikzpicture}[
    node/.style={circle, draw, minimum size=7mm, font=\footnotesize},
    edge/.style={->, >=Stealth, thick, shorten >=1pt, shorten <=1pt}
]
\node[node] (0) at (0,2) {0};
\node[node] (1) at (2,2) {1};
\node[node] (2) at (4,2) {2};
\node[node] (3) at (0,0) {3};
\node[node] (4) at (2,0) {4};
\node[node] (5) at (4,0) {5};
\node[node] (6) at (0,-2) {6};
\node[node] (7) at (2,-2) {7};
\node[node] (8) at (4,-2) {8};

\draw[edge] (0) to (2); \draw[edge] (0) to (5); \draw[edge] (0) to (8);
\draw[edge] (1) to (0); \draw[edge] (1) to (6); \draw[edge] (1) to (8);
\draw[edge] (2) to (3); \draw[edge] (2) to (4); \draw[edge] (2) to (7);
\draw[edge] (3) to (0); \draw[edge] (3) to (6); \draw[edge] (3) to (7);
\draw[edge] (4) to (1); \draw[edge] (4) to (3); \draw[edge] (4) to (7);
\draw[edge] (5) to (1); \draw[edge] (5) to (2); \draw[edge] (5) to (4);
\draw[edge] (6) to (0); \draw[edge] (6) to (2); \draw[edge] (6) to (5);
\draw[edge] (7) to (1); \draw[edge] (7) to (5); \draw[edge] (7) to (8);
\draw[edge] (8) to (3); \draw[edge] (8) to (4); \draw[edge] (8) to (6);
\end{tikzpicture}
\caption{Fig1(1)}
\label{fig:4-1}
\end{subfigure}
\hfill
\begin{subfigure}[b]{0.45\textwidth}
\centering
\begin{tikzpicture}[
    node/.style={circle, draw, minimum size=7mm, font=\footnotesize},
    edge/.style={->, >=Stealth, thick, shorten >=1pt, shorten <=1pt}
]
\foreach \angle/\label in {0/0, 36/1, 72/2, 108/3, 144/4, 180/5, 216/6, 252/7, 288/8, 324/9} {
    \node[node] (\label) at (\angle:3) {\label};
}

\draw[edge] (0) to [bend left=15] (4);
\draw[edge] (0) to [bend right=10] (5);
\draw[edge] (0) to [bend left=20] (9);
\draw[edge] (1) to [bend right=15] (3);
\draw[edge] (1) to [bend left=10] (4);
\draw[edge] (1) to [bend left=20] (7);
\draw[edge] (2) to [bend right=20] (0);
\draw[edge] (2) to [bend left=15] (1);
\draw[edge] (2) to [bend right=10] (4);
\draw[edge] (3) to [bend left=20] (2);
\draw[edge] (3) to [bend right=15] (6);
\draw[edge] (3) to [bend left=10] (8);
\draw[edge] (4) to [bend right=20] (6);
\draw[edge] (4) to [bend left=15] (8);
\draw[edge] (4) to [bend right=10] (9);
\draw[edge] (5) to [bend left=20] (1);
\draw[edge] (5) to [bend right=15] (3);
\draw[edge] (5) to [bend left=10] (7);
\draw[edge] (6) to [bend right=20] (1);
\draw[edge] (6) to [bend left=15] (2);
\draw[edge] (6) to [bend right=10] (5);
\draw[edge] (7) to [bend left=20] (0);
\draw[edge] (7) to [bend right=15] (2);
\draw[edge] (7) to [bend left=10] (9);
\draw[edge] (8) to [bend right=20] (0);
\draw[edge] (8) to [bend left=15] (5);
\draw[edge] (8) to [bend right=10] (7);
\draw[edge] (9) to [bend left=20] (3);
\draw[edge] (9) to [bend right=15] (6);
\draw[edge] (9) to [bend left=10] (8);
\end{tikzpicture}
\caption{Fig1(2)}
\label{fig:4-2}
\end{subfigure}

\label{fig:both-graphs}
\end{figure}

For $m=9$ and $m=10$, we present the graphs $\Gamma$ in Figure 1(1) and Figure 1(2) respectively. The specific edge sets of these two graphs are as follows:

For $m=9$:
\[
\begin{aligned}
\Gamma = \mathrm{Graph}[
    & 0 \to 2,\ 0 \to 5,\ 0 \to 8, \\
    & 1 \to 0,\ 1 \to 6,\ 1 \to 8, \\
    & 2 \to 3,\ 2 \to 4,\ 2 \to 7, \\
    & 3 \to 0,\ 3 \to 6,\ 3 \to 7, \\
    & 4 \to 1,\ 4 \to 3,\ 4 \to 7, \\
    & 5 \to 1,\ 5 \to 2,\ 5 \to 4, \\
    & 6 \to 0,\ 6 \to 2,\ 6 \to 5, \\
    & 7 \to 1,\ 7 \to 5,\ 7 \to 8, \\
    & 8 \to 3,\ 8 \to 4,\ 8 \to 6
]
\end{aligned}
\]

For $m=10$:
\[
\begin{aligned}
\Gamma = \mathrm{Graph}[
    & 0 \to 4,\ 0 \to 5,\ 0 \to 9, \\
    & 1 \to 3,\ 1 \to 4,\ 1 \to 7, \\
    & 2 \to 0,\ 2 \to 1,\ 2 \to 4, \\
    & 3 \to 2,\ 3 \to 6,\ 3 \to 8, \\
    & 4 \to 6,\ 4 \to 8,\ 4 \to 9, \\
    & 5 \to 1,\ 5 \to 3,\ 5 \to 7, \\
    & 6 \to 1,\ 6 \to 2,\ 6 \to 5, \\
    & 7 \to 0,\ 7 \to 2,\ 7 \to 9, \\
    & 8 \to 0,\ 8 \to 5,\ 8 \to 7, \\
    & 9 \to 3,\ 9 \to 6,\ 9 \to 8
]
\end{aligned}
\]

Clearly, $\Gamma$ is a 3-regular oriented graph in each case. Using Mathematica, we verify that both graphs are antisymmetric, i.e., their automorphism groups are trivial.
\begin{figure}[H]
  \centering
  \includegraphics[width=0.5\linewidth]{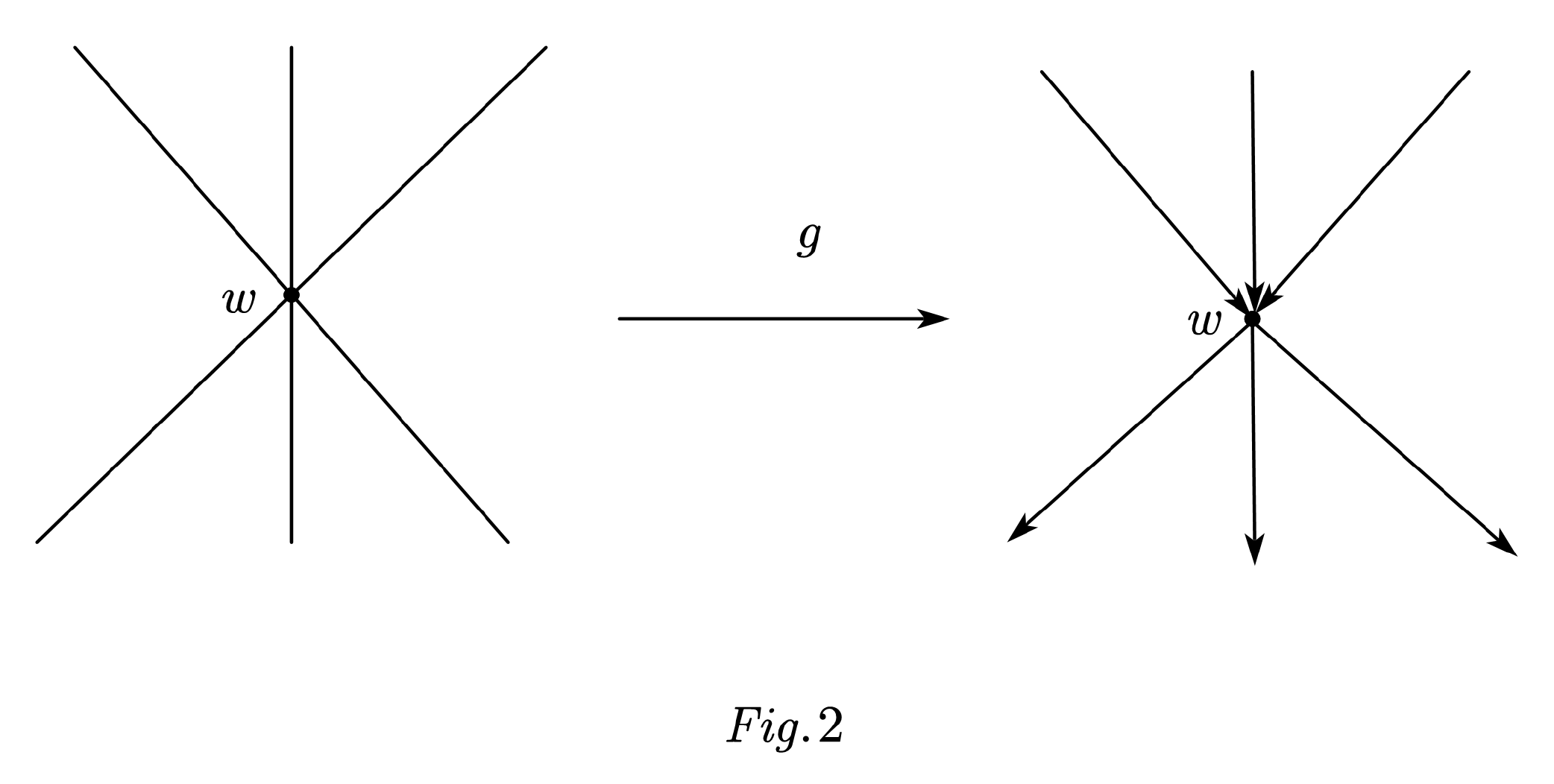}
\end{figure}

For $m \geq 11$, it follows from \cite{BG} that there exists an $m$-vertex 6-regular antisymmetric graph $\Gamma$. We perform the operation $g$ illustrated in Fig. 2: for each vertex $w$ in $\Gamma$ with its six incident edges, we convert three of them into out-edges and the remaining three into in-edges, resulting in a new graph $\Gamma' = g(\Gamma)$. Clearly, $\Gamma'$ is a 3-regular oriented graph. Since $\mathrm{Aut}(\Gamma') \subseteq \mathrm{Aut}(\Gamma) = 1$, we conclude that $\mathrm{Aut}(\Gamma') = 1$, which completes the proof.
\end{proof}

For \( G = \langle x \rangle \) with \( o(x) \geq 3 \), we discuss the following four cases: \( m = 2 \), \( m = 3 \), \( m = 4 \), and \( m \geq 5 \), which correspond to Lemma 3.3, Lemma 3.4, Lemma 3.5, and Lemma 3.6 respectively.

\begin{Lemma}
Let $G=\langle x\rangle$ be a finite group with $o(x)\geq3$.
Then $G$ admits an  2-POSR of valency 3 if and only if $o(x)\geq 7$.
\end{Lemma}
\begin{proof}
Let \( \Gamma = \mathrm{Cay}(G, T_{i,j} : i,j \in \mathbb{Z}_2) \) be an oriented 2-partite Cayley digraph with valency 3 and \( A = \mathrm{Aut}(\Gamma) \). By the definition of \( \Gamma \), we clearly have \( T_{i,i} = \emptyset \) and \( T_{i,j} \cap T_{j,i}^{-1} = \emptyset \), which will be implicitly assumed in subsequent discussions.

\textbf{Case 1.} \( o(x) \leq 6 \)

When $o(x)\leq 5$, since \( |T_{1,0}| = |T_{0,1}| = 3 \), it necessarily follows that \( T_{1,0} \cap T_{0,1}^{-1} \neq \emptyset \), which contradicts the requirement that \( \Gamma \) is an oriented graph.

When $o(x)=6$, we get $G$ does not admit a 2-POSR of valency 3 using Mathematica or MAGMA.

\textbf{Case 2.} \( o(x) = 7 \)

We define \( T_{0,1} = \{1, x, x^2\} \) and \( T_{1,0} = \{x, x^3, x^4\} \). Using Mathematicas or MAGMA\cite{MAG}, we determine that \( A = R(G) \cong G \).

\textbf{Case 3.} \( o(x) \geq 8 \)

We choose \( T_{0,1} = \{1, x, x^2\} \) and \( T_{1,0} = \{x, x^2, x^4\} \). It is evident that \( T_{0,1} \cap T_{1,0}^{-1} = \emptyset \).

For \( o(x) = 8 \) or \( 9 \), using Mathematica or MAGMA, we have \( A = R(G) \cong G \).

For \( o(x) \geq 10 \), we have:
\[
\Gamma^{+}(1_0) = \{1_1, x_1, x^2_1\}, \quad \Gamma^{+}(1_1) = \{x_0, x^2_0, x^4_0\};
\]
\[
\Gamma^{+2}(1_0) =\{x_0,x^2_0,x^4_0\}\cup \{x^2_0,x^3_0,x^5_0\}\cup \{x^3_0,x^4_0,x^6_0\}= \{x_0, x^2_0, x^3_0, x^4_0, x^5_0, x^6_0\}, \]
\[
\Gamma^{+2}(1_1) = \{x_1,x^2_1,x^3_1\}\cup \{x^2_1,x^3_1,x^4_1\}\cup \{x^4_1,x^5_1,x^6_1\}= \{x_1, x^2_1, x^3_1, x^4_1, x^5_1, x^6_1\};
\]

Next, we consider \( \Gamma^{+3}(1_i) \cap \Gamma^{+}(1_i) \) for \( i \in \mathbb{Z}_2 \).
\begin{figure}[H]
  \centering
  \includegraphics[width=0.6\linewidth]{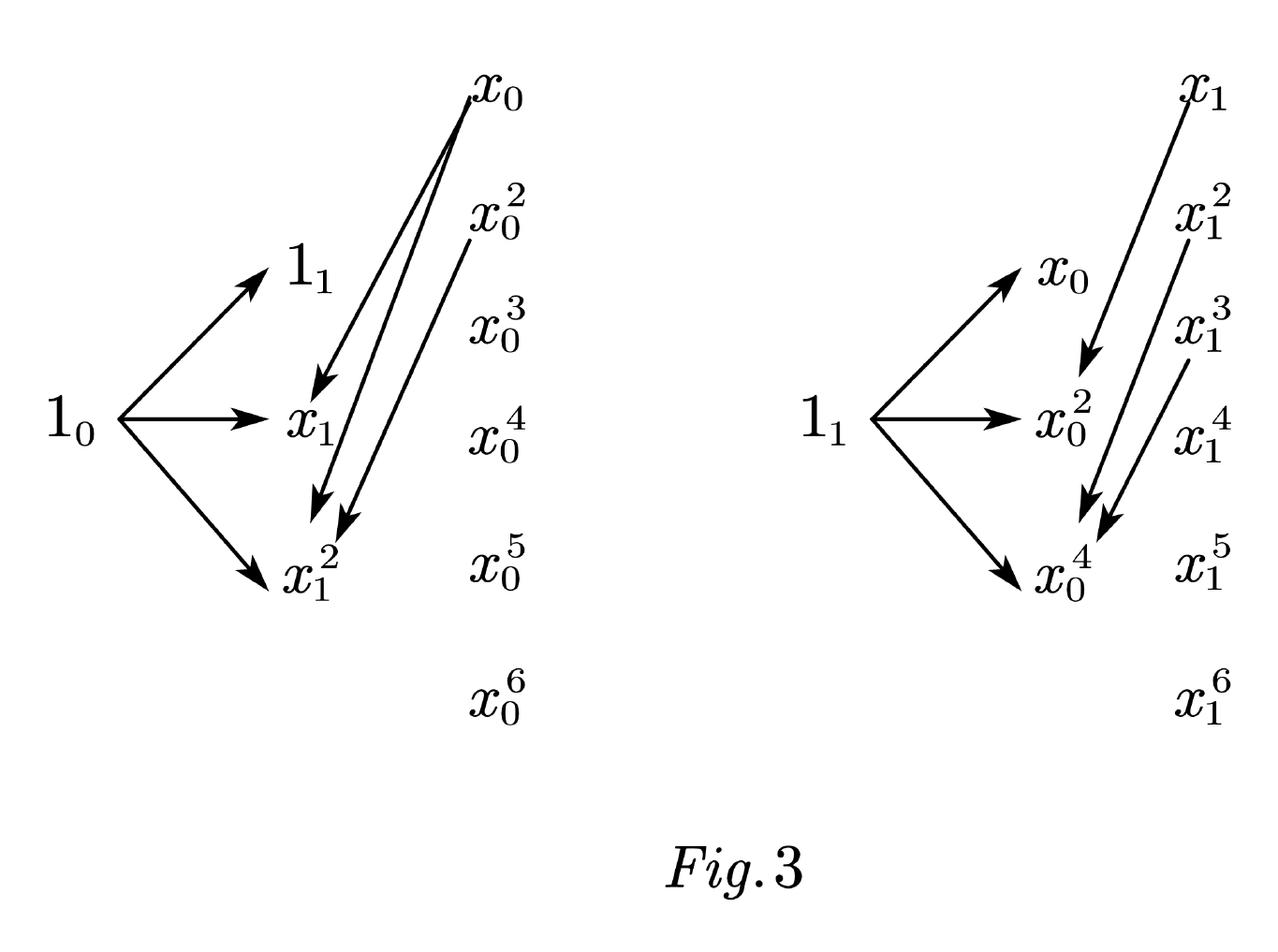}
\end{figure}
We observe (as can be directly seen in Fig. 3):
\[
\Gamma^{+}(x_0) \cap \Gamma^{+}(1_0) = \{x_1, x^2_1\}, \quad \Gamma^{+}(x^2_0) \cap \Gamma^{+}(1_0) = \{x^2_1\},
\]
\[
\Gamma^{+}(x^3_0) \cap \Gamma^{+}(1_0) = \Gamma^{+}(x^4_0) \cap \Gamma^{+}(1_0) = \Gamma^{+}(x^5_0) \cap \Gamma^{+}(1_0) = \Gamma^{+}(x^6_0) \cap \Gamma^{+}(1_0) = \emptyset;
\]
\[
\Gamma^{+}(x_1) \cap \Gamma^{+}(1_1) = \{x^2_0\}, \quad \Gamma^{+}(x^2_1) \cap \Gamma^{+}(1_1) = \{x^4_0\},
\]
\[
\Gamma^{+}(x^3_1) \cap \Gamma^{+}(1_1) = \{x^4_0\}, \quad \Gamma^{+}(x^4_1) \cap \Gamma^{+}(1_1) = \Gamma^{+}(x^5_1) \cap \Gamma^{+}(1_1) = \Gamma^{+}(x^6_1) \cap \Gamma^{+}(1_1) = \emptyset.
\]

We observe that in \( \Gamma^{+2}(1_0) \), there are two vertices \( x_0 \) and \( x^2_0 \) with out-degrees in \( \Gamma^{+}(1_0) \), whereas in \( \Gamma^{+2}(1_1) \), there are three vertices \( x_1 \), \( x^2_1 \), and \( x^3_1 \) with out-degrees in \( \Gamma^{+}(1_0) \). Therefore, \( A \) fixes \( G_0 \) and \( G_1 \) setwise.

By the Frattini argument, we have \( A = R(G)A_{g_i} \), which means:
\begin{align}
|A_{g_i}| = |A_{h_j}| = |A/R(G)|,  \tag{1}
\end{align}
for all \( g,h \in G \) and \( i,j \in \mathbb{Z}_2 \).

Additionally, the in-degrees of \( 1_1 \), \( x_1 \), and \( x^2_1 \) in \( \Gamma^{+2}(1_0) \) are 0, 1, and 2, respectively. Thus, \( A_{1_0} \) fixes every point in \( \Gamma^{+}(1_0) = \{1_1, x_1, x^2_1\} \). 
By the same method, we can show that $A_{1_1}$ fixes $\Gamma^{+}(1_1)$ pointwise, and hence $A = R(G)$ by Proposition 1.

\end{proof}

\begin{Lemma}
Let $G=\langle x\rangle$ be a finite group with $o(x)\geq3$.
Then $G$ admits an  3-POSR of valency 3 if and only if $o(x)\geq 4$.
\end{Lemma}

\begin{proof}
Let \( \Gamma = \mathrm{Cay}(G_i, T_{i,j} : i,j \in \mathbb{Z}_3) \) be an oriented 3-partite Cayley digraph with valency 3 and \( A = \mathrm{Aut}(\Gamma) \). 

\textbf{Case 1:} \( o(x) = 3 \)

When $o(x)=3$, we get $G$ does not a admit 3-POSR of valency 3 using Mathematica or MAGMA.

\textbf{Case 2:} \( o(x) \geq 4 \)

Let \( T_{0,1} = \{1, x\} \), \( T_{0,2} = \{1\} \); \( T_{1,0} = \{x^2\} \), \( T_{1,2} = \{1, x\} \); \( T_{2,0} = \{x, x^2\} \), \( T_{2,1} = \{x\} \), with all other \( T_{i,j} = \emptyset \). As shown in Fig.4, \( \Gamma \) is clearly an  oriented 3-partite Cayley digraph of valency 3. We now prove that \( A = R(G) \).
\begin{figure}[H]
  \centering
  \includegraphics[width=1.0\linewidth]{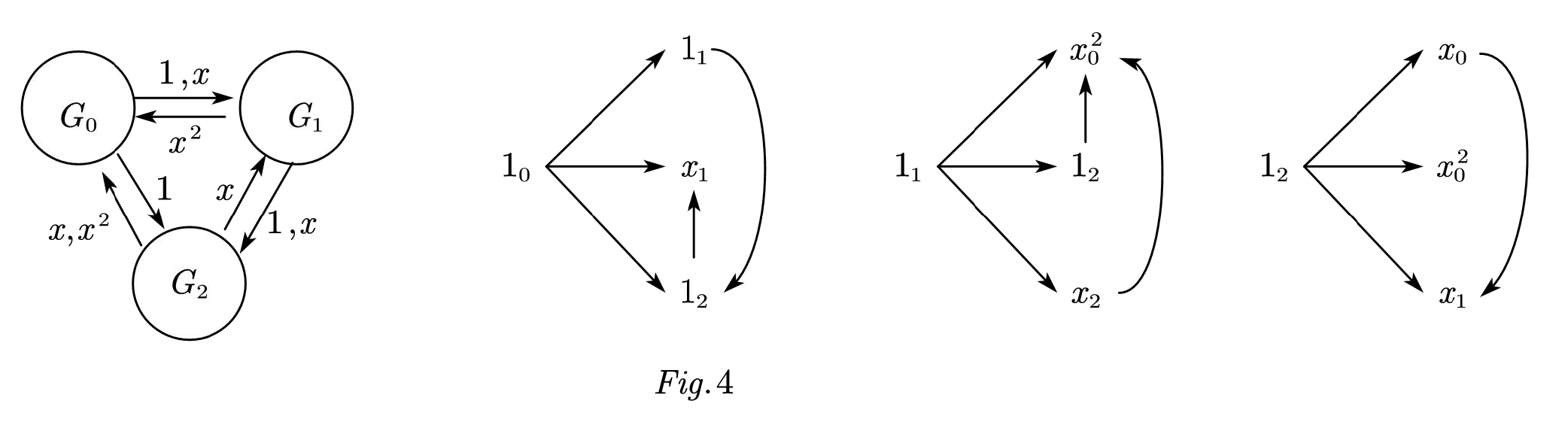}
\end{figure}

We have:
\[
\Gamma^{+}(1_0) = \{1_1, x_1, 1_2\}, \quad \Gamma^{+}(1_1) = \{x^2_0, 1_2, x_2\}, \quad \Gamma^{+}(1_2) = \{x_0, x^2_0, x_1\};
\]
\[
\Gamma^{+2}(1_0) = \Gamma^{+}(1_1) \cup \Gamma^{+}(x_1) \cup \Gamma^{+}(1_2)
= \{x^2_0, \underline{1_2}, x_2\} \cup \{x^3_0, x_2, x^2_2\} \cup \{x_0, x^2_0, \underline{x_1}\};
\]
\[
\Gamma^{+2}(1_1) = \Gamma^{+}(x^2_0) \cup \Gamma^{+}(1_2) \cup \Gamma^{+}(x_2)
= \{x^2_1, x^3_1, x^2_2\} \cup \{x_0, \underline{x^2_0}, x_1\} \cup \{\underline{x^2_0}, x^3_0, x^2_1\};
\]
\[
\Gamma^{+2}(1_2) = \Gamma^{+}(x_0) \cup \Gamma^{+}(x^2_0) \cup \Gamma^{+}(x_1)
= \{\underline{x_1}, x^2_1, x_2\} \cup \{x^2_1, x^3_1, x^2_2\} \cup \{x^3_0, x_2, x^2_2\};
\]

The subdigraphs \( [\Gamma^{+}(1_i)] \) are shown in Fig.4 for $i\in \mathbb{Z}_3$. Note that there are no isolated vertices in \( [\Gamma^{+}(1_0)] \) and \( [\Gamma^{+}(1_1)] \), but there is an isolated vertex \( x^2_0 \) in \( [\Gamma^{+}(1_2)] \). Thus, \( A \) fixes \( G_2 \) and \( G_0 \cup G_1 \) setwise. Furthermore, since vertices in \( G_0 \) and \( G_1 \) have out-degrees 1 and 2 in \( G_2 \) respectively, \( A \) fixes each \( G_i \) setwise for \( i \in \mathbb{Z}_3 \).

We next prove that \( A_{1_i} \) fixes \( \Gamma^{+}(1_i) \) pointwise for \( i \in \mathbb{Z}_3 \). Since \( A \) fixes each \( G_i \) setwise, \( A_{1_0} \) fixes \( \{1_1, x_1\} \) and \( \{1_2\} \) setwise, which means \( A_{1_0} \) fixes \( 1_2 \). As shown in Fig.4, there is an arc \( (1_1, 1_2) \) but no arc \( (x_1, 1_2) \), so \( 1_1 \) and \( x_1 \) cannot be swapped. Hence, \( A_{1_0} \) fixes \( \Gamma^{+}(1_0) = \{1_1, x_1, 1_2\} \) pointwise. Using the same reasoning, we can show that  \( A_{1_2} \) fixes \( \Gamma^{+}(1_2) \) pointwise. 
Since $A_{1_0}$ fix $\Gamma^{+}(1_0)$ pointwise, $A_{1_0}=A_{1_1}=A_{1_2}$ by Frattini argument. So $A_{1_1}$ fix $\Gamma^{+}(1_1)$ pointwise by $A$ fix $G_i$ setwise and $A_{1_1}=A_{1_2}$ for $i\in \mathbb{Z}_3$.
By Proposition 1, we conclude that \( A = R(G) \). This completes the proof.
\end{proof}

\begin{Lemma}
Let $G=\langle x\rangle$ be a finite group with $o(x)\geq3$.
Then $G$ admits an  4-POSR of valency 3.
\end{Lemma}

\begin{proof}
Define the following connection sets:
\begin{align*}
T_{0,1} &= T_{1,2} = T_{2,0} = \{1, x\}, \\
T_{3,0} &= T_{3,1} = T_{2,3} = \{1\}, \\
T_{0,3} &= T_{1,3} = \{x\}, \\
T_{3,2} &= \{x^2\}.
\end{align*}
with all other $T_{i,j} = \emptyset$ for $i,j \in \mathbb{Z}_4$. 
\begin{figure}[H]
  \centering
  \includegraphics[width=0.4\linewidth]{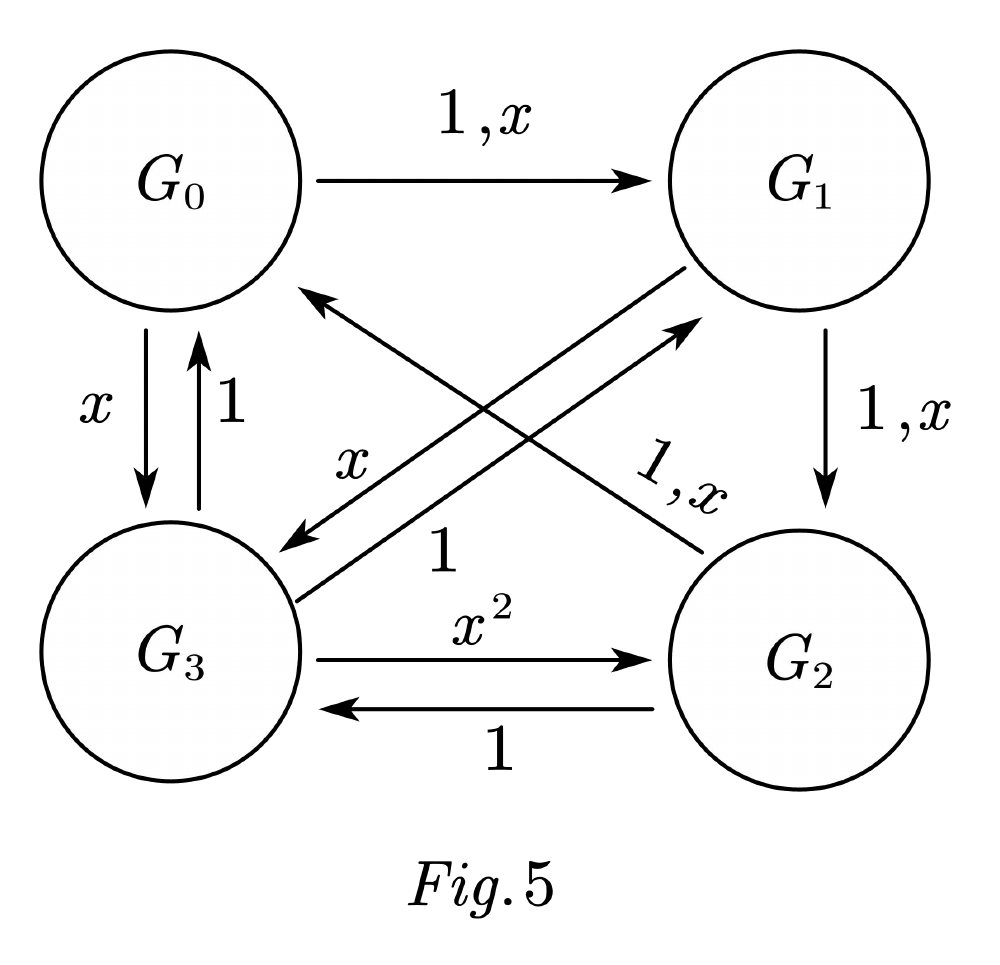}
\end{figure}
The resulting digraph $\Gamma = \mathrm{Cay}(G, T_{i,j}: i,j\in\mathbb{Z}_4)$, illustrated in Fig.5, is clearly an oriented 4-partite Cayley digraph of valency 3.
 We now prove that \( A = R(G) \), where $A=Aut(\Gamma)$.

From Fig.5, we easily observe the following: vertices in $G_i$ have their out-neighbors lying in two different parts for $i\in \{0,1,2\}$, while out-neighbors of $G_3$ lie in three different parts. Hence, $A$ fixes $G_3$ setwise. 

Furthermore, we note that $(g_3, g_1, g_2)$ forms an oriented cycle, whereas $(g_3, g_0, g_2)$ and $(g_3, g_1, g_0)$ are not oriented cycles for any $g \in G$. Thus, $A$ fixes $G_0$ setwise. Additionally, since vertices in $G_1$ and $G_2$ have out-degrees 0 and 2 respectively in $G_0$, $A$ fixes $G_1$ and $G_2$ setwise. In summary, $A$ fixes each $G_i$ setwise for $i\in \mathbb{Z}_4$.

Next, we prove that $A_{1_i}$ fixes $\Gamma^+(1_i)$ pointwise for $i\in \mathbb{Z}_4$. Since $\Gamma^+(1_3)$ lies in distinct parts, and $A$ fixes each $G_i$ setwise, $A_{1_3}$ fixes $\Gamma^+(1_3)$ pointwise. Because $T_{3,0} = T_{3,1} = T_{2,3} = \{1\}$, we have $A_{1_3} = A_{1_0} = A_{1_1} = A_{1_2}$. 

Given $\Gamma^+(1_0) = \{1_1, x_1, x_3\}$, and since $A$ fixes each $G_i$ setwise with $A_{1_0} = A_{1_1}$, it follows that $A_{1_0}$ fixes $\Gamma^+(1_0)$ pointwise.

For $\Gamma^+(1_1) = \{x_3, 1_2, x_2\}$ and $\Gamma^+(1_2) = \{1_3, 1_0, x_0\}$, using the fact that $A$ fixes each $G_i$ setwise and $A_{1_3} = A_{1_0} = A_{1_1} = A_{1_2}$, we similarly conclude that $A_{1_1}$ fixes $\Gamma^+(1_1)$ pointwise and $A_{1_2}$ fixes $\Gamma^+(1_2)$ pointwise.

In summary, $A_{1_i}$ fixes $\Gamma^+(1_i)$ pointwise for all $i \in \mathbb{Z}_3$. By Proposition 1, we obtain $A = R(G)$.

\end{proof}

\begin{Lemma}
Let $G=\langle x\rangle$ be a finite group with $o(x)\geq3$ .
Then $G$ admits an  $m$-POSR of valency 3 for $m\geq5$.
\end{Lemma}
\begin{proof}
We define the connection sets $T_{i,j} \subseteq G$ as follows:
\begin{align*}
T_{i,i+1} &= \{1\}, \quad T_{i,i-1} = \{x\} \quad \text{for } i \in \mathbb{Z}_m \\
T_{j,j-2} &= \{1\}, \quad T_{2,0} = \{x\} \quad \text{for } j \in \mathbb{Z}_m\backslash\{2\}
\end{align*}
with all other $T_{i,j} = \emptyset$ for $i,j \in \mathbb{Z}_m$.
Note that here we assume \( o(x) \geq 3 \) since \( G \ncong \mathbb{Z}_2 \).
\begin{figure}[H]
  \centering
  \includegraphics[width=1.0\linewidth]{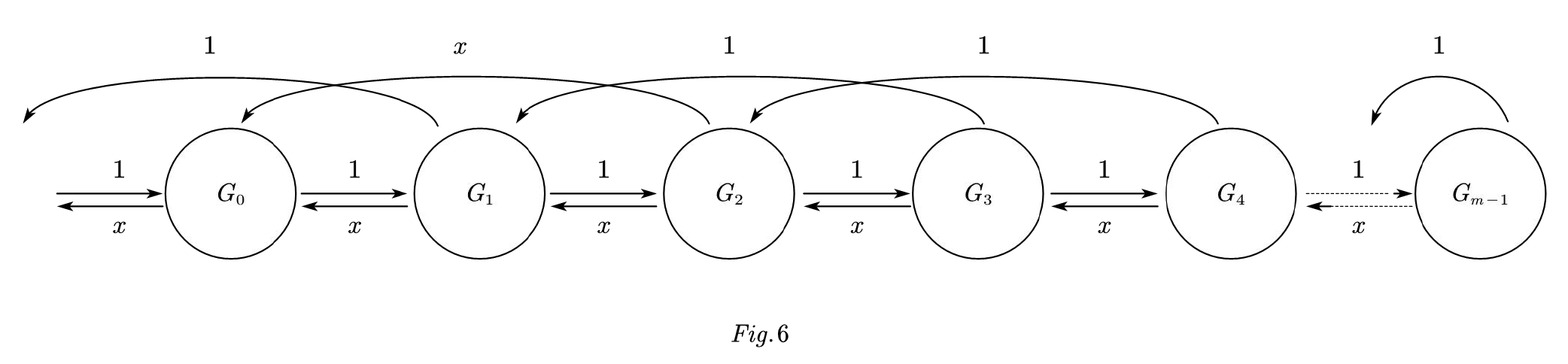}
\end{figure}
The resulting digraph $\Gamma = \mathrm{Cay}(G, T_{i,j}:i,j\in\mathbb{Z}_m)$, as shown in Fig.6, is clearly an oriented $m$-partite Cayley digraph of valency 3. Let $\mathrm{Aut}(\Gamma) = A$; we now prove that $A = R(G)$.

From Fig.6, we observe that $(g_i, g_{i+1}, g_{i+2})$ forms a oriented 3-cycle for all $i \neq 0$. This implies that $A$ fixes $G_0 \cup G_1 \cup G_2$ setwise. Examining the induced subgraph $H = [G_0 \cup G_1 \cup G_2]$, we note that:

\begin{itemize}
\item Vertices in $G_0$ and $G_1$ have in-degree 2 in $H$.
\item Vertices in $G_2$ have in-degree 1 in $H$.
\end{itemize}
Consequently, $A$ fixes $G_2$ setwise. Moreover, since $T_{1,2} = \{1\}$ and $T_{0,2} = \emptyset$, $A$ must also fix $G_0$ and $G_1$ setwise.

Since \( T_{2,i} = \emptyset \) for \( i \notin \{0,1,3\} \), \( A \) fixes \( G_3 \) setwise because \( A \) fixes \( G_0, G_1, G_2 \) setwise. Similarly, since \( T_{3,i} = \emptyset \) for \( i \notin \{1,2,4\} \), \( A \) fixes \( G_4 \) setwise because \( A \) fixes \( G_1, G_2, G_3 \) setwise. By induction, we conclude that \( A \) fixes every \( G_i \) setwise for \( i \in \mathbb{Z}_m \).

Note that the vertices in \( \Gamma^+(1_i) \) lie in distinct parts. Therefore, since \( A \) fixes each \( G_i \) setwise, \( A \) fixes \( \Gamma^+(1_i) \) pointwise for all \( i \in \mathbb{Z}_m \). By Proposition 1, we thus have \( A = R(G) \).

\end{proof}

\begin{proof}[Proof of Theorem 1.2]
Let \( \Gamma = \mathrm{Cay}(G, T_{i,j} : i,j \in \mathbb{Z}_2) \) be an oriented 2-partite Cayley digraph with valency 3 and \( A = \mathrm{Aut}(\Gamma) \). 
\subsection*{Case 1: When $o(y)=2$}
We have $o(xy)=2$. If $o(xy)\geq 3$, then we discuss in Case 2 by $G=\langle x,y\rangle=\langle x,xy\rangle$. So $xyx=y^{-1}$, i.e., 
\[
G=\langle x,y \mid o(x)=4, o(y)=2, xyx=y^{-1}\rangle\cong D_8.
\]
Let $T_{0,1}=\{1,x,xy\}$, $T_{1,0}=\{x,y,x^3y\}$, then we have $A=R(G)\cong G$ by Mathematica or MAGMA.

\subsection*{Case 2: When $o(y)\geq 3$, i.e., $o(y)=3$ or $4$}
We take $T_{0,1}=\{1,x,y\}$, $T_{1,0}=\{x,x^2,y\}$.
\begin{align*}
\Gamma^{+}(1_0) &= \{1_1,x_1,y_1\}, \\
\Gamma^{+}(1_1) &= \{x_0,x^2_0,y_0\}, \\
\Gamma^{+2}(1_0) &= \Gamma^{+}(1_1)\cup \Gamma^{+}(x_1)\cup \Gamma^{+}(y_1) \\
&= \{x_0,x^2_0,y_0\}\cup \{x^2_0,x^3_0,(yx)_0\}\cup \{(xy)_0,(x^2y)_0,y^2_0\}, \\
\Gamma^{+2}(1_1) &= \Gamma^{+}(x_0)\cup \Gamma^{+}(x^2_0)\cup \Gamma^{+}(y_0) \\
&= \{x_1,x^2_1,(yx)_1\}\cup \{x^2_1,x^3_1,(yx^2)_1\}\cup \{y_1,(xy)_1,y^2_1\}.
\end{align*}
To find $\Sigma=[\Gamma^{+}[1_0]\cup \Gamma^{+2}(1_0)]$ and $\Delta=[\Gamma^{+}[1_1]\cup \Gamma^{+2}(1_1)]$, we still need to find $\Gamma^{+3}(1_0)$ and $\Gamma^{+3}(1_1)$, as follows:
\begin{align*}
\Gamma^{+}(x^3_0) &= \{x^3_1,1_1,(yx^3)_1\}, \\
\Gamma^{+}((yx)_0) &= \{(yx)_1,\underline{(xyx)_1},(y^2x)_1\}, \\
\Gamma^{+}((xy)_0) &= \{(xy)_1,(x^2y)_1,\underline{(yxy)_1}\}, \\
\Gamma^{+}((x^2y)_0) &= \{(x^2y)_1,(x^3y)_1,\underline{(yx^2y)_1}\}, \\
\Gamma^{+}(y^2_0) &= \{y^2_1,(xy^2)_1,\underline{y^3_1}\}, \\
\Gamma^{+}((yx)_1) &= \{\underline{(xyx)_0},\underline{(x^2yx)_0},(y^2x)_0\}, \\
\Gamma^{+}(x^3_1) &= \{1_0,x_0,(yx)_0\}, \\
\Gamma^{+}((yx^2)_1) &= \{\underline{(xyx^2)_0},\underline{(x^2yx^2)_0},(y^2x^2)_0\}, \\
\Gamma^{+}((xy)_1) &= \{(x^2y)_0,(x^3y)_0,\underline{(yxy)_0}\}.
\end{align*}
It is worth noting that the underlined elements may lead to the following situations:
\begin{enumerate}
    \item $\Gamma^{+}((yx)_0)\cap \Gamma^{+}(1_0)=(xyx)_1=y_1$
    \item $\Gamma^{+}((xy)_0)\cap \Gamma^{+}(1_0)=(yxy)_1=x_1$
    \item $\Gamma^{+}((x^2y)_0)\cap \Gamma^{+}(1_0)=(yx^2y)_1=1_1$ or $x_1$
    \item $\Gamma^{+}(y^2_0)\cap \Gamma^{+}(1_0)=y^3_1=1_1$
    \item $\Gamma^{+}((yx)_1)\cap \Gamma^{+}(1_1)=(xyx)_0=y_0$
    \item $\Gamma^{+}((yx)_1)\cap \Gamma^{+}(1_1)=(x^2yx)_0=y_0$
    \item $\Gamma^{+}((yx^2)_1)\cap \Gamma^{+}(1_1)=(xyx^2)_0=y_0$
    \item $\Gamma^{+}((yx^2)_1)\cap \Gamma^{+}(1_1)=(x^2yx^2)_0=y_0$
    \item $\Gamma^{+}((xy)_1)\cap \Gamma^{+}(1_1)=(yxy)_0=x_0$ or $x^2_0$
\end{enumerate}
We write all the above special cases as a set:
\[
H=\{xyx=y, yxy=x, yx^2y=1, yx^2y=x, x^2yx=y, x^2yx^2=y, xyx^2=y, yxy=x^2, y^3=1\}
\]
\[
=\{y=xyx, y=xy^{-1}x^3, y=y^{-1}x^2, y=xy^{-1}x^2, y=x^2yx, y=x^2yx^2, y=xyx^2, y=x^2y^{-1}x^3, y^3=1\}.
\]
If $y=x^2yx$ holds, then $yx^{-1}y^{-1}=x^2$, i.e., $|x^{-1}|=|x^2|$, which contradicts $o(x)=4$. For $y = x y x^2$, we can similarly derive a contradiction. In what follows, we let $H = \{y=xyx, y=xy^{-1}x^3, y=y^{-1}x^2, y=xy^{-1}x^2, y=x^2yx^2, y=x^2y^{-1}x^3, y^3=1\}$.

\begin{figure}[H]
  \centering
  \includegraphics[width=1.0\linewidth]{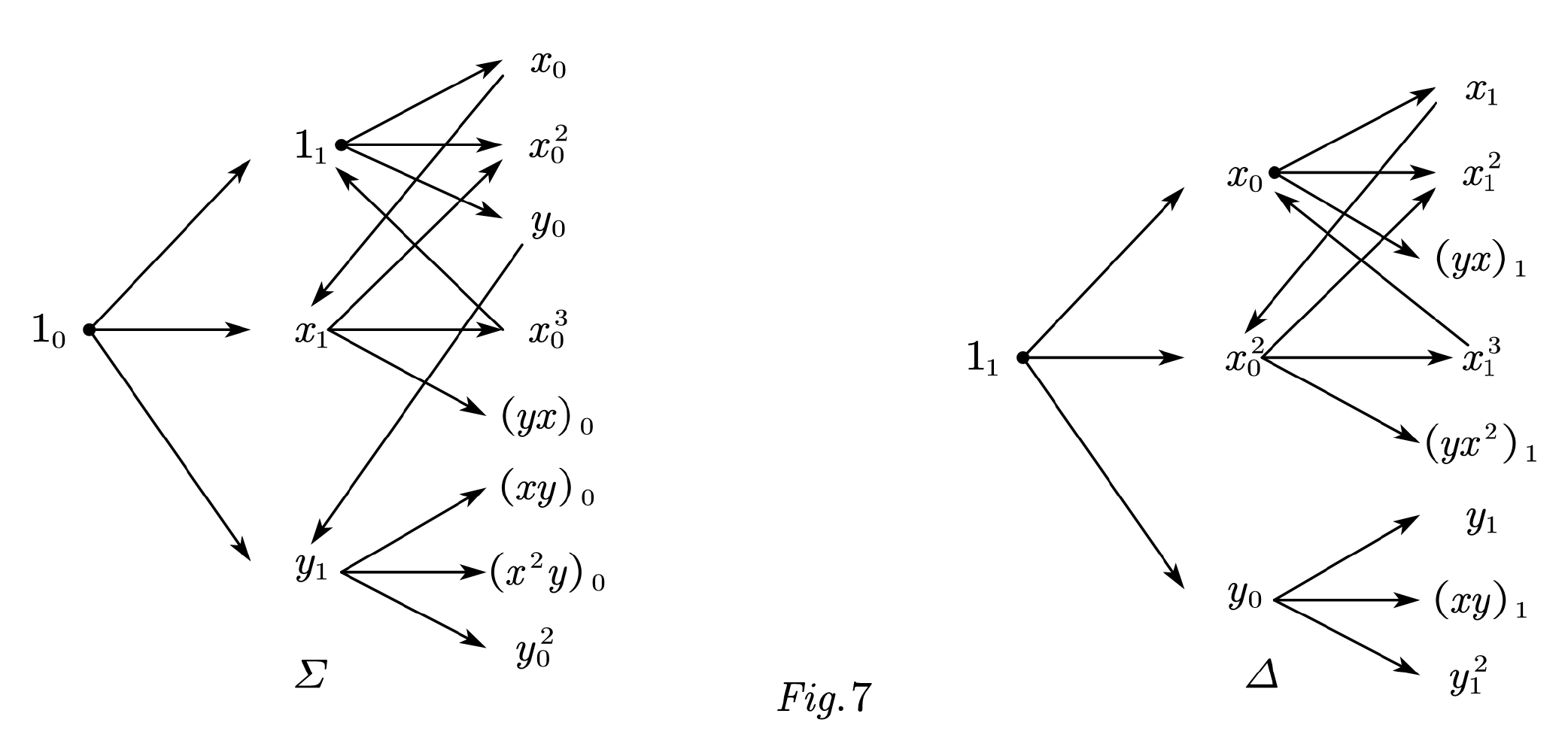}
\end{figure}

\begin{figure}[H]
  \centering
  \includegraphics[width=1.0\linewidth]{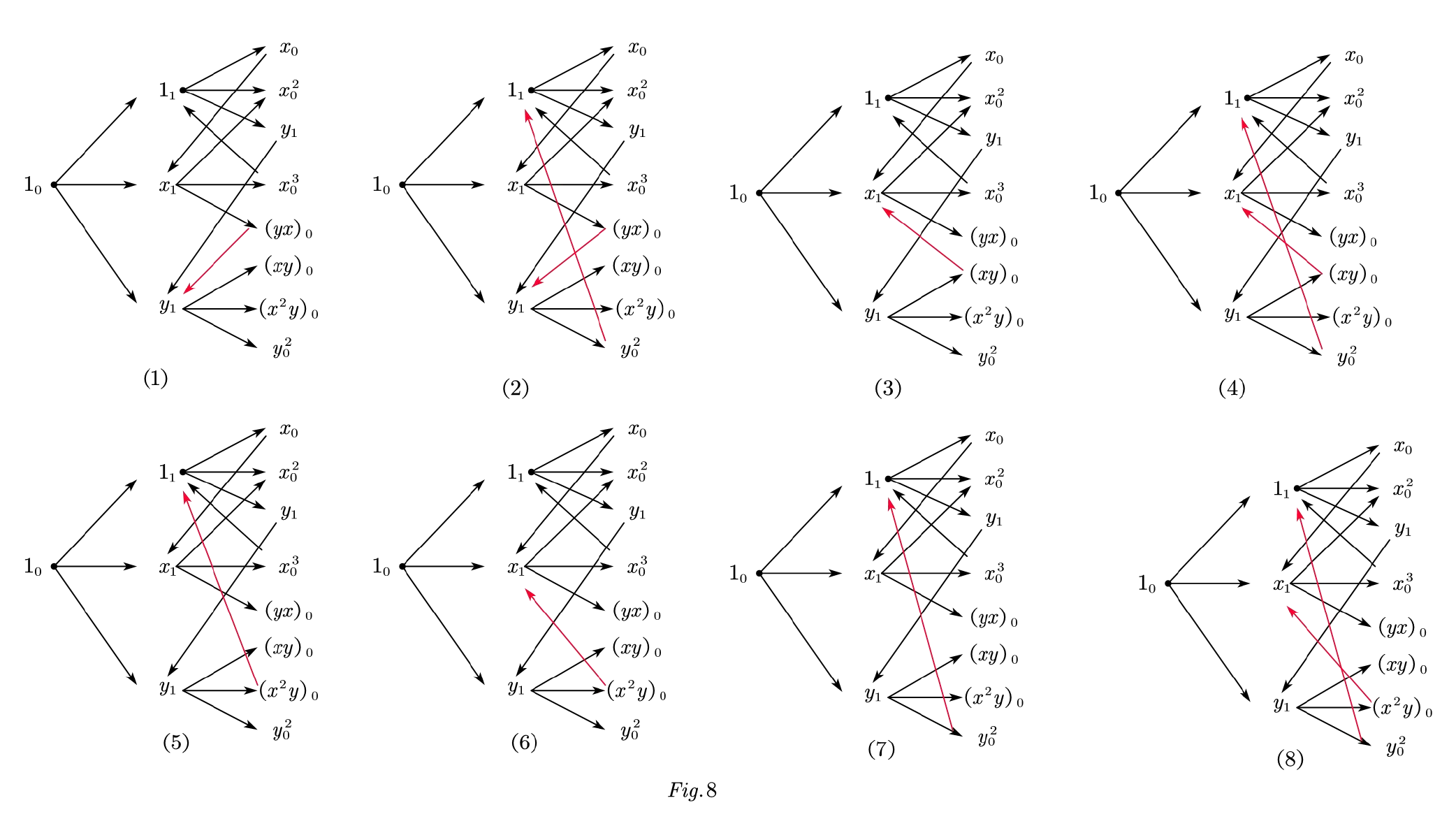}
\end{figure}

\begin{figure}[H]
  \centering
  \includegraphics[width=1.0\linewidth]{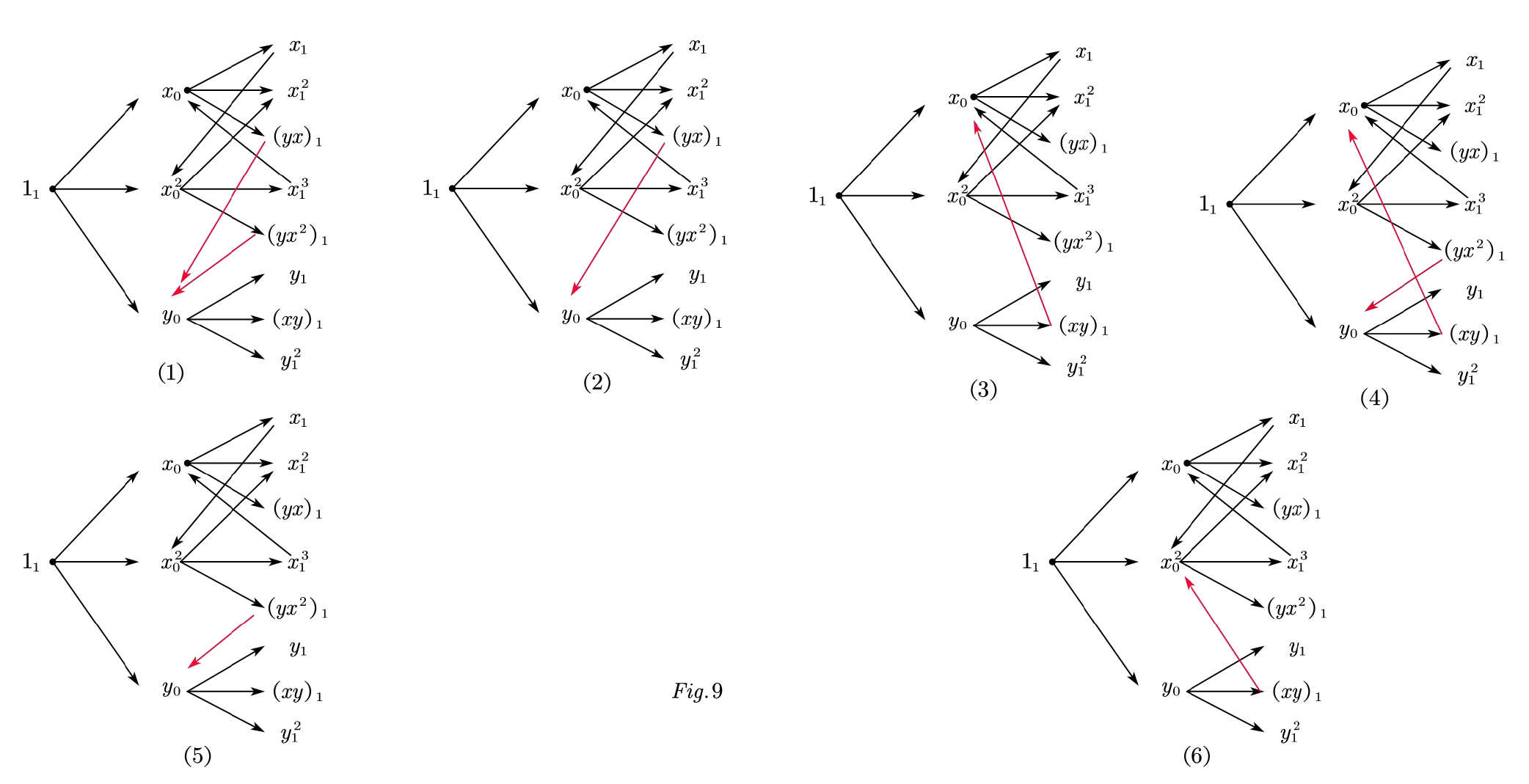}
\end{figure}

\subsubsection*{Case 2.1: If none of the equations in $H$ hold}
Then $\Sigma$ and $\Delta$ are as shown in Fig.8. Clearly $\Sigma\ncong \Delta$, so we have $A$ fixes $G_0$, $G_1$ setwise. $A_{1_0}$ clearly fixes $\{1_1,x_1\}$ and $\{y_1\}$ setwise.$(1_0,1_1,y_0,y_1)$ is the unique oriented $P_4$ in $\Sigma$ passing through $y_1$, so $A_{1_0}$ fixes $1_1$, $y_0$, and $y_1$ since $A_{1_0}$ fixes $y_1$. So $A_{1_0}$ fixes $\Gamma^{+}(1_0)=\{1_1,x_1,y_1\}$ pointwise. By the Frattini argument we know $A_{1_0}=A_{1_1}=A_{x_1}=A_{y_1}$, so $A_{1_0}$ fixes $\Gamma^{+}(1_1)=\{x_0,x^2_0,y_0\}$ setwise. The unique in-neighbor of $x_1$ in $\Gamma^{+2}(1_0)$ is $x_0$, so $A_{1_0}$ fixes $x_0$ by $A_{1_0}=A_{x_1}$.
Consequently, $A_{1_1}$ fixes $x_0$ because $A_{1_0} = A_{1_1}$.
From $\Delta$ we clearly see $A_{1_1}$ fixes $y_0$, so we have $A_{1_1}$ fixes $\Gamma^{+}(1_1)$ pointwise. In summary, we have $A=R(G)$ by Proposition 1.

Now we start discussing the cases where some equations in $H$ hold.

\subsubsection*{Case 2.2: $y=xyx$}
\begin{enumerate}
    \item $y=xy^{-1}x^3$: We have $xyx=xy^{-1}x^3$ i.e., $y^2=x^2$, so $|y|=|x|=4$, then $G=Q_8$. Using Mathematica or MAGMA we find that $Q_8$ does not have a 2-POSR of valency 3.
    \item $y\neq xy^{-1}x^3$, $y=y^{-1}x^2$: We have $y^2=x^2$ i.e., $|y|=|x|=4$, then again $G=Q_8$.
    \item $y\neq xy^{-1}x^3$, $y\neq y^{-1}x^2$, $y=xy^{-1}x^2$: Then $y=xyx=xy^{-1}x^2$, so $y^2=x$, contradiction.

    \item $y\neq xy^{-1}x^3$, $y\neq y^{-1}x^2$, $y\neq xy^{-1}x^2$,  $y=x^2y^{-1}x^3$: Then $xyx = x^2y^{-1}x^3$, which implies $y = xy^{-1}x^2 = xyx$, so $x = y^2$, a contradiction.
    \item $y\neq xy^{-1}x^3$, $y\neq y^{-1}x^2$, $y\neq xy^{-1}x^2$,   $y\neq x^2y^{-1}x^3$, $y=x^2yx^2$: Then we have $\Delta$ as in Fig.9(1). If $y^3 \neq 1$, i.e., $y^4 = 1$, then we have $G = C_4\rtimes C_4$. If $y^3 = 1$, then $\Sigma$ is as shown in Fig.~8(2), and we easily observe that $\Sigma \ncong \Delta$, so $A$ fixes $G_0$ and $G_1$ setwise. In this case, it is straightforward to verify that $A_{1_0}$ fixes $\{1_1, x_1\}$ and $\{y_1\}$ setwise. Since the in-degrees of $1_1$ and $x_1$ in $\Gamma^{+2}(1_0)$ are $2$ and $1$ respectively, $A_{1_0}$ fixes $\Gamma^{+}(1_0)$ pointwise, i.e., $A_{1_0} = A_{1_1} = A_{x_1} = A_{y_1}$. In $\Sigma$, $x_0$ is the unique in-neighbor of $x_1$ in $\Gamma^{+2}(1_0)$, so $A_{1_0}$ fixes $x_0$ because $A_{1_0} = A_{x_1}$. Hence, $A_{1_1}$ fixes $x_0$ as $A_{1_0} = A_{1_1}$. Moreover, $A_{1_1}$ clearly fixes $y_0$, so $A_{1_1}$ fixes $\Gamma^{+}(1_1) = \{x_0, x^2_0, y_0\}$ pointwise. Therefore, $A = R(G)$ by Proposition 1.
    \item $y\neq xy^{-1}x^3$, $y\neq y^{-1}x^2$, $y\neq xy^{-1}x^2$,  $y\neq x^2y^{-1}x^3$, $y\neq x^2yx^2$, $y^3=1$: Then $\Sigma$ as in Fig.8(2), $\Delta$ as in Fig.9(2). Similar to the proof in Case 2.2.5, we can get $A=R(G)$.
\end{enumerate}
This completes the discussion for $y=xyx$. Below we assume $y\neq xyx$.

\subsubsection*{Case 2.3:$y = x y^{-1} x^3$ (equivalent to $y x y = x$)}
\begin{enumerate}
    \item $y=y^{-1}x^2$: Then we have $|x|=|y|=4$, then $G=Q_8$ by $yxy=x$ and $y^2=x^2$.
    \item $y\neq y^{-1}x^2$, $y=xy^{-1}x^2$: Then $xy^{-1}x^3=xy^{-1}x^2$ so $x=1$, contradiction.
    \item $y\neq y^{-1}x^2$, $y\neq xy^{-1}x^2$,$y=x^2y^{-1}x^3$: Then $xy^{-1}x^3=x^2y^{-1}x^3$ so $x=1$, contradiction.
    
     Now we still have $y=x^2yx^2$ or $y^3=1$ not discussed in $H$. Now we discuss them:

        \item $y\neq y^{-1}x^2$, $y\neq xy^{-1}$, $y\neq xyx^2$, $y\neq x^2y^{-1}x^3$,  $y=x^2yx^2$: Then $\Delta$ as in Fig.9(4). If also $y^3=1$, then $\Sigma$ as in Fig.8(4), similar to the proof in Case 2.2.5, we have $A=R(G)$.
        If $y^3 \neq 1$, then $y^4 = 1$, and we have $G = C_4\rtimes C_4$.
        \item $y\neq y^{-1}x^2$, $y\neq xy^{-1}$, $y\neq xyx^2$, $y\neq x^2y^{-1}x^3$,  $y\neq x^2yx^2$, $y^3=1$: Then $\Sigma$ as in Fig. 8(4), $\Delta$ as in Fig.9(3). Similarly we can prove $A=R(G)$.
        \item $y\neq y^{-1}x^2$, $y\neq xy^{-1}$, $y\neq xyx^2$, $y\neq x^2y^{-1}x^3$,  $y\neq x^2yx^2$, $y^3\neq 1$: Then $\Sigma$ as in Fig.8(3), $\Delta$ as in Fig.9(3). Similarly we have $A=R(G)$.
   
\end{enumerate}
This completes the discussion for $y\neq xyx$, $y=xy^{-1}x^3$. Below we assume $y\neq xyx$ and $y\neq xy^{-1}x^3$.

\subsubsection*{Case 2.4: $y=y^{-1}x^2$ (equivalent to $y x^2 y = 1$)}
Then clearly $o(y)=4$ by $x^2=y^2$.
\begin{enumerate}
    \item $y=xy^{-1}x^2$: Then $y^{-1}x^2=xy^{-1}x^2$ so $x=1$, contradiction.
    \item $y\neq xy^{-1}x^2$,  $y=x^2y^{-1}x^3$: Then $y=x^2y^{-1}x^3=y^2y^{-1}x^3$ so $x^3=1$, contradiction.
    \item $y\neq xy^{-1}x^2$, $y\neq x^2y^{-1}x^3$, $y^3=1$: This is clearly a contradiction because $o(y)=4$.
    \item $y\neq xy^{-1}x^2$, $y\neq x^2y^{-1}x^3$, $y^3\neq1$, $y=x^2yx^2$: Then $\Sigma$ as in Fig.8(5), $\Delta$ as in Fig.9(5). Similar to Case 2.2.6, we can prove $A=R(G)$.
    \item $y\neq xy^{-1}x^2$, $y\neq x^2y^{-1}x^3$, $y^3\neq1$, $y\neq x^2yx^2$: Then $\Sigma$ as in Fig.8(5), $\Delta$ as in Fig.7. Similarly we have $A=R(G)$.
\end{enumerate}

Now we start considering $y\neq xyx$, $y\neq xy^{-1}x^3$ and $y\neq y^{-1}x^2$.

\subsubsection*{Case 2.5: $y=xy^{-1}x^2$ (equivalent to $y x^2 y = x$)}
\begin{enumerate}
    \item $y=x^2yx^2$: Then $xy^{-1}x^2=x^2yx^2$ so $x=y^{-2}$, contradiction.
    \item $y\neq x^2yx^2$, $y=x^2y^{-1}x^3$: Then $xy^{-1}x^2=x^2y^{-1}x^3$, so $y^{-1}=xy^{-1}x$, so $y=xy^{-1}x^2=(xy^{-1}x)x=y^{-1}x$, so $x=y^2$, contradiction.
    \item $y\neq x^2yx^2$, $y\neq x^2y^{-1}x^3$, $y^3=1$: Then $\Sigma$ as in Fig.8(8), $\Delta$ as in Fig.7. Similarly we can get $A=R(G)$.
    \item $y\neq x^2yx^2$, $y\neq x^2y^{-1}x^3$, $y^3\neq 1$: Then $\Sigma$ as in Fig.8(6), $\Delta$ as in Fig.7. Similarly we can get $A=R(G)$.
\end{enumerate}

Now we start considering $y\neq xyx$, $y\neq xy^{-1}x^3$, $y\neq y^{-1}x^2$, $y\neq xy^{-1}x^2$.

\subsubsection*{Case 2.6: $y=x^2yx^2$}
\begin{enumerate}
    \item $y=x^2y^{-1}x^3$: Then $x^2y^{-1}x^3=x^2yx^2$ so $x=y^2$, contradiction.
    \item  $y\neq x^2y^{-1}x^3$, $y^3=1$: Then $\Sigma$ as in Fig.8(7), $\Delta$ as in Fig.9(5). Similarly we have $A=R(G)$.
    \item $y\neq x^2y^{-1}x^3$, $y^3\neq 1$: Then $y^4=1$, we defer the discussion of this case to the end.  
\end{enumerate}

Now we start considering $y\neq xyx$, $y\neq xy^{-1}x^3$, $y\neq y^{-1}x^2$, $y\neq xy^{-1}x^2$, $y\neq x^2yx^2$.

\subsubsection*{Case 2.8: $y=x^2y^{-1}x^3$(equivalent to $y x y = x^2$)}
\begin{enumerate}
    \item $y^3=1$: Then $\Sigma$ as in Fig.8(7), $\Delta$ as in Fig.9(6). Similarly we have $A=R(G)$.
    \item $y^3\neq 1$: Then $\Sigma$ as in Fig.7, $\Delta$ as in Fig.9(6). Similarly we have $A=R(G)$.
\end{enumerate}

Now we start considering $y\neq xyx$, $y\neq xy^{-1}x^3$, $y\neq y^{-1}x^2$, $y\neq xy^{-1}x^2$, $y\neq x^2yx^2$, $y\neq x^2y^{-1}x^3$.

\subsubsection*{Case 2.9: $y^3=1$}
Then $\Sigma$ as in Fig.8(7), $\Delta$ as in Fig.7. Similarly we have $A=R(G)$.

At this point, only one case remains undiscussed, namely  
\[
G = \langle x, y \mid o(x) = o(y) = 4,\ y = x^2 y x^2 \rangle.
\]  
Since $G = \langle x, xy \rangle$, we have $o(xy) \neq 3$; otherwise, we return to Case 2.6.2, which has already been discussed. Hence, $o(xy) \in \{2, 4\}$. Similarly, $o(yx) \in \{2, 4\}$ and $o(x^2 y) \in \{2, 4\}$.

Define a mapping $f$ such that $f(x) = y$, $f(y) = x$, $f(T_{0,1}) = T'_{0,1}$, and $f(T_{1,0}) = T'_{0,1}$. This yields a new graph $\Gamma' = \mathrm{Cay}(G, T'_{i,j},\ i,j \in \mathbb{Z}_2)$. By symmetry, we have $\Gamma \cong \Gamma'$, so $\Gamma'$ also exhibits the same special case:  
\[
G = \langle x, y \mid o(x) = o(y) = 4,\ x = y^2 x y^2 \rangle.
\]  
Therefore, the case we need to examine is:  
\[
G = \langle x, y \mid o(x) = o(y) = 4,\ y = x^2 y x^2,\ x = y^2 x y^2 \rangle.
\]

If $o(xy) = 2$ or $o(yx) = 2$, then  
\[
G = \langle x, y \mid o(x) = o(y) = 4,\ o(xy) = 2,\ y = x^2 y x^2,\ x = y^2 x y^2 \rangle \cong \mathrm{SmallGroup}(16, 3) = (C_2 \times C_2) \rtimes C_4.
\]

Now assume $o(xy) = o(yx) = 4$, i.e.,  
\[
G = \langle x, y \mid o(x) = o(y) = o(xy) = o(yx) = 4,\ y = x^2 y x^2,\ x = y^2 x y^2 \rangle.
\]  
Note that this could still be an infinite group, so we further consider $o(x^2 y)$.  

If $o(x^2 y) = 2$, then from $y = x^2 y x^2$ we have $x^2 y = y x^2$. Left-multiplying both sides by $x^2 y$ gives  
\[
(x^2 y)(x^2 y) = (x^2 y)(y x^2) \Rightarrow 1 = x^2 y^2 x^2,
\]  
so $y^2 = 1$, contradicting $y^4 = 1$. Hence, $o(x^2 y) = 4$.  

In this case, we obtain  
$G = \langle x, y \mid o(x) = o(y) = o(xy) = o(yx) = o(x^2 y) = 4,\ y = x^2 y x^2,\ x = y^2 x y^2 \rangle \cong \mathrm{SmallGroup}(32, 2) = (C_4 \times C_2) \rtimes C_4$.
For $G\in \{Q_8, C_4\rtimes C_4, (C_2\times C_2)\rtimes C_4, (C_4\times C_2)\rtimes C_4$\}, we know G does not admit a 2-POSR of valency 3 using Mathematica or MAGMA.

So far we have discussed all cases.

\end{proof}

\begin{proof}[Proof of Theorem 1.3]
When the order of any generator in $G$ is less than or equal to 3, we know that $G \in \{ \mathbb{Z}_2^2, \mathbb{Z}_3^2, D_6, A_4, He_3 \}$ by \cite[4.1]{du4} where $D_n$ denotes the dihedral group of order $n$, $He_3$ denotes the Heisenberg group of order 27, and $A_n$ denotes the alternating group of degree $n$.. 

For $G = \mathbb{Z}_3^2 = \langle x, y \mid o(x) = o(y) = 3,\ xy = yx \rangle$, take $T_{0,1} = \{1, x, y\}$ and $T_{1,0} = \{y, x, xy^2\}$. Then, using Mathematica or MAGMA, we obtain $A = R(G)$.

For $G = A_4 = \langle x, y \mid x^3 = 1, y^2 = 1, (xy)^3 = 1 \rangle$, take $T_{0,1} = \{1, yx, yxy\}$ and $T_{1,0} = \{1, yx, xyxy\}$. Then, using Mathematica or MAGMA, we obtain $A = R(G)$.

For $G = \mathrm{He}_3 = \langle x, y \mid x^3 = y^3 = z^3 = 1, [x,y] = z, [x,z] = [y,z] = 1 \rangle$, take $T_{0,1} = \{z, x^2z^2, x^2yz\}$ and $T_{1,0} = \{x^2z^2, z, x^2y^2\}$. Then, using Mathematica or MAGMA, we obtain $A = R(G)$.
Using Mathematica or MAGMA, we know that in this case $G$ does not admit a 2-POSR of valency 3 when $G=D_6$ or $\mathbb{Z}^2_2$.

Next, we discuss the case where $G$ has a generator $a$ with order $\geq 5$. Without loss of generality, assume $o(x) \geq 5$.

\noindent \textbf{Case 1:} If $G$ does not have a generator of order $\geq 6$.

\noindent \textbf{Case 1.1:} $o(x) = 5$, $o(y) = 2$.

We have $G = \langle x, y \rangle = \langle x, xy \rangle$, so $o(xy) = 2$; otherwise $G$ would fall into Case 1.2. Therefore, $xyx = x^{-1}$, i.e., $G \cong D_{10}$. 

Take $T_{0,1} = \{1, x, x^2\}$ and $T_{1,0} = \{x, y, xy\}$. Then, using Mathematica or MAGMA, we obtain $A = R(G) \cong G$..

\noindent \textbf{Case 1.2:} $o(x) = 5$, $o(y) \geq 3$.

Take $T_{0,1} = \{1, x, x^2\}$, $T_{1,0} = \{x, y, y^2\}$.
\begin{figure}[H]
  \centering
  \includegraphics[width=1.0\linewidth]{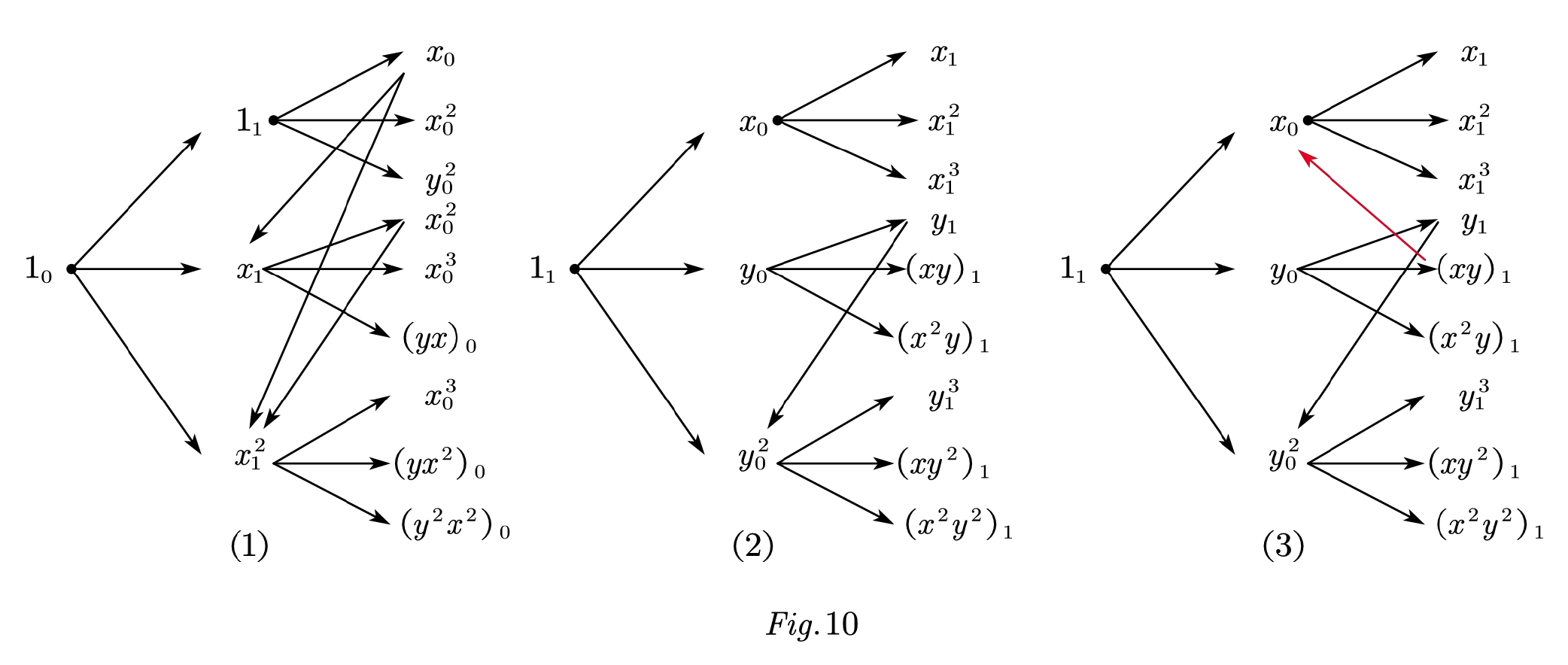}
\end{figure}
Then:
\begin{align*}
\Gamma^{+}(1_0) &= \{1_1, x_1, x^2_1\}, \\
\Gamma^{+}(1_1) &= \{x_0, y_0, y^2_0\}, \\
\Gamma^{+2}(1_0) &= \Gamma^{+}(1_1) \cup \Gamma^{+}(x_1) \cup \Gamma^{+}(x^2_1) \\
&= \{x_0, y_0, y^2_0\} \cup \{x^2_0, (yx)_0, (y^2x)_0\}\cup \{ x^3_0, yx^2, y^2x^2\}, \\
\Gamma^{+2}(1_1) &= \Gamma^{+}(x_0) \cup \Gamma^{+}(y_0) \cup \Gamma^{+}(y^2_0) \\
&= \{x_1, x^2_1, x^3_1\} \cup \{y_1, (xy)_1, (x^2y)_1\} \cup \{y^2_1, (xy^2)_1, (x^2y^2)_1\}.
\end{align*}

To find $\Sigma = [\Gamma^{+}[1_0] \cup \Gamma^{+2}(1_0)]$ and $\Delta = [\Gamma^{+}[1_1] \cup \Gamma^{+2}(1_1)]$, we also need to compute $\Gamma^{+3}(1_0)$ and $\Gamma^{+3}(1_1)$. However, to avoid excessive detail, we omit the computation of $\Gamma^{+3}(1_0)$ and $\Gamma^{+3}(1_1)$ and directly present some special cases:

\begin{align*}
\Gamma^{+}((xy)_1) &= \{(x^2y)_0, \underline{(yxy)_0}, \underline{(y^2xy)_0}\}, \\
\Gamma^{+}((x^2y)_1) &= \{(x^3y)_0, \underline{(yx^2y)_0}, \underline{(y^2x^2y)_0}\}, \\
\Gamma^{+}((xy^2)_1) &= \{(x^2y^2)_0, \underline{(yxy^2)_0}, \underline{(y^2xy^2)_0}\}, \\
\Gamma^{+}((x^2y^2)_1) &= \{(x^3y^2)_0, \underline{(yx^2y^2)_0}, \underline{(y^2x^2y^2)_0}\}.
\end{align*}

Note that $x$ might belong to the set 
\[
H = \{yxy, y^2xy, yx^2y, y^2x^2y, yxy^2, y^2xy^2, yx^2y^2, y^2x^2y^2\}.
\]
We can easily prove that it is impossible for two elements in $H$ to both equal $x$.

If no element in $H$ equals $x$, then $\Sigma$ and $\Delta$ are as shown in Fig. 10(1) and Fig. 10(2). When an element in $H$ equals $x$, without loss of generality, assume $yxy = x$. Then $\Delta$ is as shown in Fig.10(3).

From Fig.10, it is easy to see that regardless of whether $x$ is in $H$, we have $\Sigma \ncong \Delta$, so $A$ fixes $G_0$ and $G_1$ setwise.

Since the in-degrees of $1_1$, $x_1$, $x^2_1$ in $\Gamma^{+2}(1_0)$ are 1, 1, and 2 respectively, $A_{1_0}$ fixes $\{1_1, x_1\}$ and $\{x^2_1\}$ setwise. The in-degrees of $1_1$ and $x_1$ in $\Gamma^{+}(x^2_1)$ are 1 and 0 respectively, so $A_{1_0}$ fixes $\Gamma^{+}(1_0)$ pointwise.

If $x \notin H$, then $(1_1, y_0, y_1, y^2_0)$ is the only directed path of length 3 in $\Delta$, so $A_{1_1}$ fixes $y_0$ and $y^2_0$, i.e., $A_{1_1}$ fixes $\Gamma^{+}(1_1)$ pointwise.

When $yxy = x$, the in-degrees of $x_0$, $y_0$, $x^2_0$ in $\Gamma^{+2}(1_1)$ are 1, 0, and 1 respectively, so $A_{1_1}$ fixes $y_0$. Since $A_{1_0} = A_{1_1}$, and in $\Sigma$, the in-degrees of $x_0$, $y_0$, $x^2_0$ in $\Gamma^{+}(1_0)$ are 2, 0, and 0 respectively, $A_{1_1}$ fixes $x_0$, i.e., $A_{1_1}$ fixes $\Gamma^{+}(1_1)$ pointwise. Therefore, by Proposition 2.1, we know that $A = R(G) \cong G$.

\noindent \textbf{Case 2: $6\leq o(x)$} 
\begin{figure}[H]
  \centering
  \includegraphics[width=1.0\linewidth]{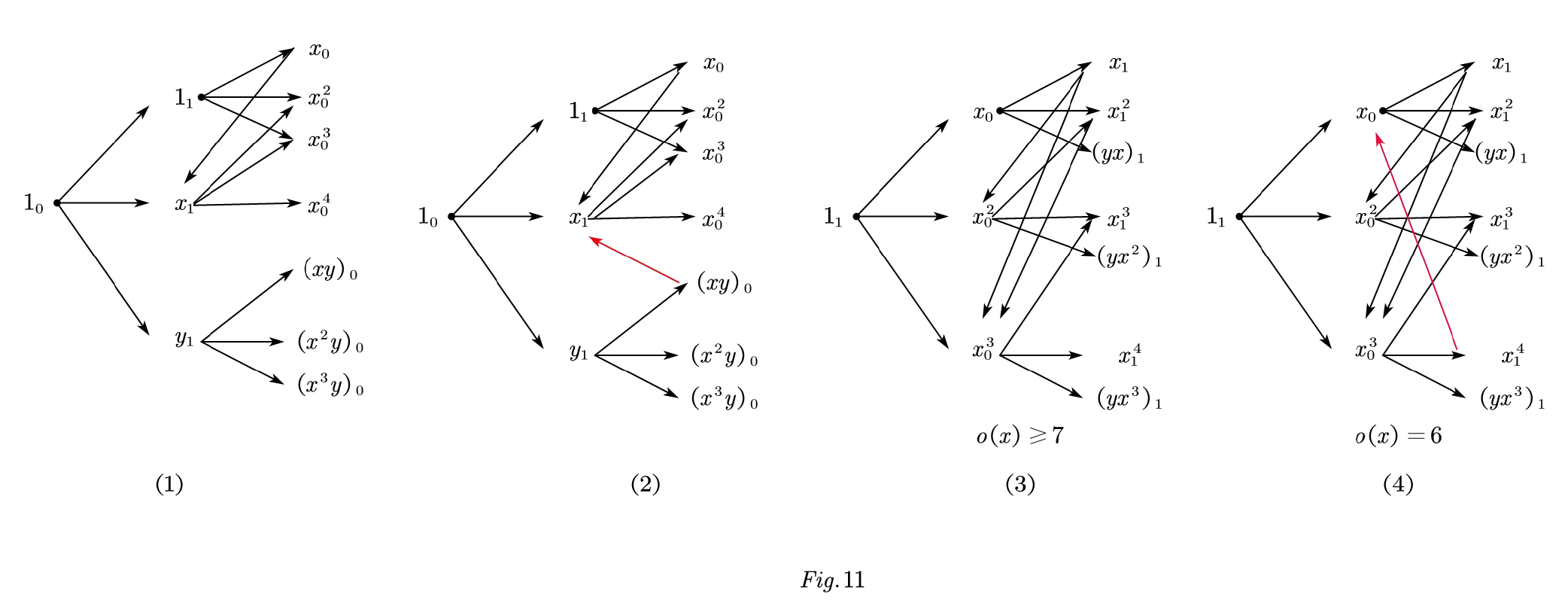}
\end{figure}
When $o(x) \geq 6$, take $T_{0,1} = \{1, x, y\}$ and $T_{1,0} = \{x, x^2, x^3\}$. Let $\Sigma = [\Gamma^{+}[1_0] \cup \Gamma^{+2}(1_0)]$ and $\Delta = [\Gamma^{+}[1_1] \cup \Gamma^{+2}(1_1)]$.

When $o(x) \geq 7$, $\Delta$ is as shown in Fig.11(3); when $o(x) = 6$, $\Delta$ is as shown in Fig.11(4). Note that
\[
\Gamma^{+}((xy)_0) = \{(xy)_1, (x^2y)_1, \underline{(yxy)_1}\}, \]
\[
\Gamma^{+}((x^2y)_0) = \{(x^2y)_1, (x^3y)_1, \underline{(yx^2y)_1}\}, \]
\[
\Gamma^{+}((x^3y)_0) = \{(x^3y)_1, (x^4y)_1, \underline{(yx^3y)_1}\},
\]
so it is possible that $x = yxy$, $x = yx^2y$, or $x = yx^3y$ holds, but not all simultaneously. If $x \notin \{yxy, yx^2y, yx^3y\}$, then $\Sigma$ is as shown in Fig.~11(1); if $x \in \{yxy, yx^2y, yx^3y\}$, without loss of generality, assume $x = yxy$, then $\Sigma$ is as shown in Fig.11(2).

When $o(x) \geq 7$, from Fig.11 we see that regardless of whether $x = xyx$ or $x \neq xyx$, we have $\Sigma \ncong \Delta$, so $A$ fixes $G_1$ and $G_2$ setwise. Note that, as shown in Fig.~11(3), the in-degrees of $x_0$, $x^2_0$, and $x^3_0$ in $\Gamma^{+2}(1_1)$ are $0$, $1$, and $2$, respectively, so $A_{1_1}$ fixes $\Gamma^{+}(1_1) = \{x_0, x^2_0, x^3_0\}$ pointwise. Also note that, as shown in Fig.~11(1) and Fig.~11(2), regardless of whether $x \neq xyx$ or $x = xyx$, we have $A_{1_0}$ fixes $\{1_1, x_1\}$ and $\{y_1\}$ setwise. Furthermore, in Fig.~11(1) or Fig.~11(2), the out-degree of $\Gamma^{+}(1_1)$ in $\Gamma^{+}(1_0)$ is $1$ (since $(x_0, x_1)$ is an arc), while the out-degree of $\Gamma^{+}(x_1)$ in $\Gamma^{+}(1_0)$ is $0$, so $A_{1_0}$ fixes $\Gamma^{+}(1_0)$ pointwise. By Proposition~1, we have $A = R(G)$. For the case $o(x) = 6$, a similar argument shows that $A = R(G)$.
\end{proof}

\begin{proof}[Proof of Theorem 1.4]
Define the connection sets as follows:
\begin{align*}
T_{i,i+1} &= \{1,x\} , 
T_{m-1,0} = \{x,y\} \text{ for } i \in \mathbb{Z}_m\backslash\{m-1\} \\
T_{i,i-1} &= \{x\} \text{ for } i \in \mathbb{Z}_m,
\end{align*}
with all other $T_{i,j} = \emptyset$.
\begin{figure}[H]
  \centering
  \includegraphics[width=1.0\linewidth]{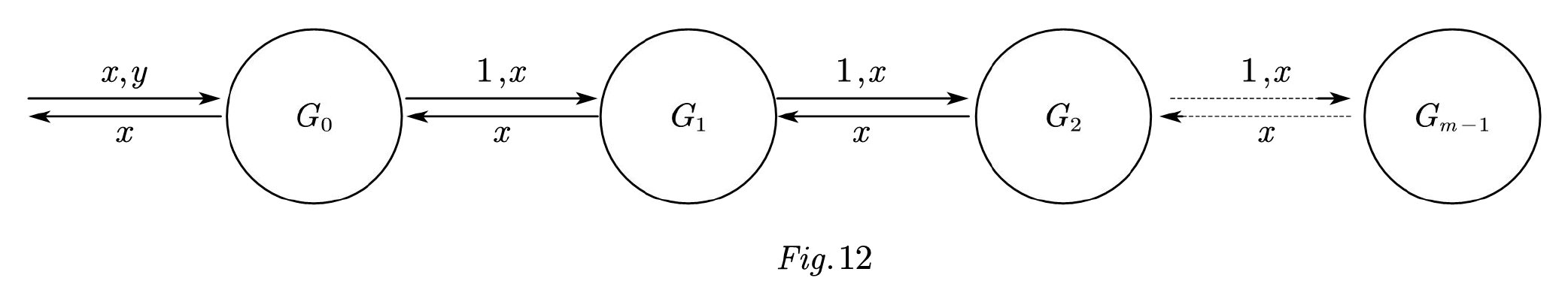}
\end{figure}
The digraph $\Gamma = \mathrm{Cay}(G, T_{i,j}:i,j\in\mathbb{Z}_m)$, shown in Fig.5, is clearly an oriented $m$-partite Cayley digraph of valency 3. We now prove $A = R(G)$.

We have:
\[
\Gamma^+(1_i) = \{1_{i+1}, x_{i+1}, x_{i-1}\} \quad \text{for } i \in \mathbb{Z}_m \setminus \{m-1\};
\]
\[
\Gamma^+(1_{m-1}) = \{x_0, y_0, x_{m-2}\};
\]

\[
\Gamma^{+2}(1_0) = \Gamma^+(1_1) \cup \Gamma^+(x_1) \cup \Gamma^+(x_{m-1})
\]
\[
= \{1_2, x_2, x_0\} \cup \{x_2, x^2_2, x^2_0\} \cup \{x^2_0, (yx)_0, x^2_{m-2}\}
\]
\[
= \{x_0, x^2_0, (yx)_0, 1_2, x_2, x^2_2, x^2_{m-2}\};
\]

\[
\Gamma^{+2}(1_i) = \Gamma^+(1_{i+1}) \cup \Gamma^+(x_{i+1}) \cup \Gamma^+(x_{i-1}) \quad \text{for } 1 \leq i \leq m-3
\]
\[
= \{1_{i+2}, x_{i+2}, x_i\} \cup \{x_{i+2}, x^2_{i+2}, x^2_i\} \cup \{1_i, x_i, x_{i-2}\}
\]
\[
= \{x_{i-2}, x_i, x^2_i, 1_i, 1_{i+2}, x_{i+2}, x^2_{i+2}\};
\]

\[
\Gamma^{+2}(1_{m-2}) = \Gamma^+(1_{m-1}) \cup \Gamma^+(x_{m-1}) \cup \Gamma^+(x_{m-3})
\]
\[
= \{x_0, y_0, x_{m-2}\} \cup \{x^2_0, (yx)_0, x^2_{m-2}\} \cup \{x_{m-2}, x^2_{m-2}, x^2_{m-4}\}
\]
\[
= \{x_0, y_0, x^2_0, (yx)_0, x_{m-2}, x^2_{m-2}, x^2_{m-4}\};
\]

\[
\Gamma^{+2}(1_{m-1}) = \Gamma^+(x_0) \cup \Gamma^+(y_0) \cup \Gamma^+(x_{m-2})
\]
\[
= \{x_1, x^2_1, x^2_{m-1}\} \cup \{y_1, (xy)_1, (xy)_{m-1}\} \cup \{x_{m-1}, x^2_{m-1}, x^2_{m-3}\}
\]
\[
= \{x_1, x^2_1, y_1, (xy)_1, x_{m-1}, x^2_{m-1}, (xy)_{m-1}, x^2_{m-3}\};
\]
Since $|\Gamma^{+2}(1_i)| = 7$ for $i \neq m-1$ and $|\Gamma^{+2}(1_{m-1})| = 8$, $A$ fixes $G_{m-1}$ setwise. The conditions $T_{m-1,i} = \emptyset$ for $i \notin \{0,m-2\}$ imply $A$ fixes $G_0 \cup G_{m-2}$ setwise. Moreover, since $|T_{m-1,0}| = 2$ and $|T_{m-1,m-2}| = 1$, $A$ fixes both $G_0$ and $G_{m-2}$ setwise.

Similarly, from $T_{i,j} = \emptyset$ for $j \notin \{i-1,i+1\}$ with $|T_{i,i+1}| = 2$ and $|T_{i,i-1}| = 1$, we conclude $A$ fixes every $G_i$ setwise.

We now show $A_{1_i}$ fixes $\Gamma^{+}(1_i)$ pointwise. Since $A$ fixes all $G_i$ setwise, $A_{1_0}$ fixes both $\{1_1,x_1\}$ and $\{x_{m-1}\}$ setwise, implying $A_{1_0} = A_{x_{m-1}}$. Thus $A_{1_0}$ fixes $\Gamma^{+}(x_{m-1})$ setwise. Observing that $\Gamma^{+}(1_1) \cap \Gamma^{+}(x_{m-1}) = \emptyset$ while $\Gamma^{+}(x_1) \cap \Gamma^{+}(x_{m-1}) = \{x^2_0\}$, vertices $1_1$ and $x_1$ cannot be swapped. Therefore $A_{1_0}$ fixes $\Gamma^{+}(1_0) = \{1_1,x_1,x_{m-1}\}$ pointwise.

By the Frattini argument, we obtain:
\begin{equation}
A_{1_0} = A_{1_1} = A_{x_1} = A_{x_{m-1}}. \tag{4}
\end{equation}

This implies:
\begin{equation}
A_{x_0} = A_{x_1} = A_{x^2_1} = A_{x^2_{m-1}}. \tag{5}
\end{equation}

Combining (4) and (5) yields:
\begin{equation}
A_{1_0} = A_{1_1} = A_{x_0} = A_{x_1}. \tag{6}
\end{equation}

From $\Gamma^{+}(1_1) = \{x_0,1_2,x_2\}$, we see $A_{1_1}$ fixes $\{x_0\}$ and $\{1_2,x_2\}$ setwise. Since $T_{2,1} = \{x\}$, we have $A_{1_2} = A_{x_1}$, which with (6) gives $A_{1_1} = A_{1_2}$. Therefore $A_{1_1}$ fixes $\Gamma^{+}(1_1)$ pointwise, establishing:
\begin{equation}
A_{1_1} = A_{x_0} = A_{1_2} = A_{x_2}. \tag{7}
\end{equation}

Combining (6) and (7) gives:
\begin{equation}
A_{1_0} = A_{1_1} = A_{1_2} = A_{x_0} = A_{x_1} = A_{x_2}. \tag{8}
\end{equation}

By induction, we generalize to:
\begin{equation}
A_{1_0} = A_{1_1} = \cdots = A_{1_{m-1}} = A_{x_0} = A_{x_1} = \cdots = A_{x_{m-1}}. \tag{9}
\end{equation}

Since $\Gamma^{+}(1_i) = \{x_{i-1},1_{i+1},x_{i+1}\}$ for $i \in \mathbb{Z}_m\backslash\{m-1\}$ and $\Gamma^{+}(1_{m-1}) = \{x_0,y_0,x_{m-2}\}$, equation (9) shows $A_{1_i}$ fixes $\Gamma^{+}(1_i)$ pointwise for all $i \in \mathbb{Z}_m$. Proposition 2.1 then establishes $A = R(G)$.
\end{proof}

\begin{proof}[Proof of Corollary 1.6]
From Theorems 1.1 and 1.5, we know that when $G=\langle x \rangle$, we only need to examine the case where $G=1$ and $m \leq 8$. For $G=\langle x,y \rangle \neq \langle x \rangle$, we only need to examine the case where $G \in \{ \mathbb{Z}_2^2, D_6, Q_8, C_4\rtimes C_4, (C_2\times C_2)\rtimes C_4, (C_4\times C_2)\rtimes C_4 \}$.

\noindent \textbf{Case 1: $G=1$ and $m \leq 8$}

For $m \leq 5$, it is clear that no $m$-vertex 3-regular digraph exists. For $m=6$, using Mathematica, we find that there is no 6-vertex 3-regular antisymmetric digraph.
For $m=7$ and $m=8$, we explicitly construct the following 3-regular digraphs:

For $m=7$:
\[
\begin{aligned}
\Gamma_7 = \mathrm{Graph}[
    & 0 \to 6,\ 0 \to 4,\ 0 \to 3, \\
    & 1 \to 4,\ 1 \to 3,\ 1 \to 6, \\
    & 2 \to 5,\ 2 \to 0,\ 2 \to 1, \\
    & 3 \to 6,\ 3 \to 4,\ 3 \to 1, \\
    & 4 \to 1,\ 4 \to 5,\ 4 \to 2, \\
    & 5 \to 2,\ 5 \to 0,\ 5 \to 3, \\
    & 6 \to 2,\ 6 \to 0,\ 6 \to 5
    ]
\end{aligned}
\]
For $m=8$:
\[
\begin{aligned}
\Gamma_8 = \mathrm{Graph}[
    & 1 \to 2,\ 1 \to 5,\ 1 \to 7, \\
    & 2 \to 3,\ 2 \to 6,\ 2 \to 8, \\
    & 3 \to 1,\ 3 \to 4,\ 3 \to 7, \\
    & 4 \to 2,\ 4 \to 5,\ 4 \to 8, \\
    & 5 \to 1,\ 5 \to 6,\ 5 \to 7, \\
    & 6 \to 3,\ 6 \to 4,\ 6 \to 8, \\
    & 7 \to 2,\ 7 \to 4,\ 7 \to 6, \\
    & 8 \to 1,\ 8 \to 3,\ 8 \to 5
   ]
\end{aligned}
\]

Clearly, both $\Gamma_7$ and $\Gamma_8$ are 3-regular digraphs. Using Mathematica, we compute $\mathrm{Aut}(\Gamma_7) = 1$ and $\mathrm{Aut}(\Gamma_8) = 1$, which implies that both digraphs are antisymmetric.

 \subsection*{Case 2: $G \in \{ \mathbb{Z}_2^2, D_6, Q_8, C_4\rtimes C_4, (C_2\times C_2)\rtimes C_4, (C_4\times C_2)\rtimes C_4 \}$}

When $G = \mathbb{Z}_2^2$, it is known from \cite{du5} that $G$ does not admit a 2-PDR.

For $D_6 = \langle x, y \mid o(x) = 3, o(y) = 2, o(xy) = 2 \rangle = \{1, x, x^2, y, xy, x^2y\}$, taking $T_{0,1} = \{1, x, x^2\}$ and $T_{1,0} = \{1, y, xy\}$, we obtain $A \cong G$ using Mathematica or MAGMA.

For $G=<x,y|o(x)=4> \in\{Q_8, C_4\rtimes C_4, (C_2\times C_2)\rtimes C_4, (C_4\times C_2)\rtimes C_4\}$, taking $T_{0,1} = \{1, x, y\}$ and $T_{1,0} = \{1, x^{-1}, x^{-2}\}$, we obtain $A \cong G$ using Mathematica or MAGMA.

\end{proof}

\section{Acknowledgments}

We gratefully acknowledge the support of the Graduate Innovation Program of China University of Mining and Technology (2025WLKXJ146), the Fundamental Research Funds for the Central Universities, and the Postgraduate Research and Practice Innovation Program of Jiangsu Province (KYCX25\_2858) for this work.


\begin{thebibliography}{s2}
\bibitem{bab1} L. Babai, Finite digraphs with given regular automorphism groups, Period. Math. Hung. 11 (1980) 257-270.
\bibitem{BG} G. Baron, W. Imrich, Asymmetrische reguläre graphen. Acta Mathematica Hungarica, 1969, 20(1-2): 135-142.
\bibitem{du1} J.-L. Du, Y.-Q. Feng, P. Spiga, A classification of the graphical $m$-semiregular representations of finite groups, J. Comb. Theory, Ser. A 171 (2020) 105174.
\bibitem{du1'} J.-L. Du, Y.-Q. Feng, S. Bang, On oriented $m$-semiregular representations of finite groups, Journal of Graph Theory, 2024, 107(3): 485-508.
\bibitem{du2} J.-L. Du, Y. Kwon, D. Yang, On oriented $m$-semiregular representations of finite groups about valency two, Discrete Math. 346 (2023) 113279.
\bibitem{du3} J.-L. Du, Y. Kwon, F. Yin, On $m$-partite oriented semiregular representations of finite groups generated by two elements, Discrete Math 347 (2024) 114043.
\bibitem{du4} J.-L. Du, Y.-Q. Feng, B. Xia, D.-W. Yang, The existence of $m$-Haar graphical representations, J. Combin. Theory Ser. A \textbf{218} (2026), Paper No. 106096, 30 pp.
\bibitem{du5} J.-L. Du, Y.Q. Feng, P. Spiga, On $n$-partite digraphical representations of finite groups. Journal of Combinatorial Theory, Series A, 2022, 189: 105606.

\bibitem{fru} R. Frucht, Graph of degree three with a given abstract group, Can. J. Math. 1 (1949) 365-378.
\bibitem{god1}C.D. Godsil, The automorphism groups of some cubic Cayley graphs, European Journal of Combinatorics, 1983, 4(1): 25-32.
\bibitem{he} D. Hetzel, \"Uber regul\"are graphische Darstellung von aufl\"osbaren Gruppen, Technische Universit\"at, Berlin, 1976 (Diplomarbeit).
\bibitem{im} W. Imrich, Graphical regular representations of groups odd order, in: Combinatorics, in: Coll. Math. Soc. J\'anos. Bolayi, vol. 18, 1976, pp. 611-621.
\bibitem{im1} W. Imrich, M.E. Watkins, On graphical regular representations of cyclic extensions of groups, Pac. J. Math. 55 (1974) 461-477.
\bibitem{im2} W. Imrich, M.E. Watkins, On automorphism groups of Cayley graphs, Period. Math. Hung. 7 (1976) 243-258.
\bibitem{kon} D. K$\ddot{o}$nig, Theory of Finite and Infinite Graphs, translated from the German by Richard McCoart, with a commentary by W. T. Tutte and a biographical sketch by T. Gallai, Birkh$\ddot{a}$user Boston, Inc., Boston, MA, 1990.
\bibitem{MAG} W. Bosma, C. Cannon, C. Playoust, The MAGMA algebra system I: the user language, J. Symb. Comput. 24 (1997) 235-265.
\bibitem{mor} J. Morris, P. Spiga, Every finite non-solvable group admits an oriented regular representation, J. Comb. Theory, Ser. B 126 (2017) 198-234.
\bibitem{mor1} J. Morris, P. Spiga, Classification of finite groups that admit an oriented regular representation, Bull. Lond. Math. Soc. 50 (2018) 811-831.
\bibitem{no} L.A. Nowitz, M.E. Watkins, Graphical regular representations of non-abelian groups I, Can. J. Math. 24 (1972) 994-1008.
\bibitem{no1} L.A. Nowitz, M.E. Watkins, Graphical regular representations of non-abelian groups II, Can. J. Math. 24 (1972) 1009-1018.
\bibitem{spi1} P. Spiga, Cubic graphical regular representations of finite non-abelian simple groups, Comm. Algebra 46 (2018), 2440-2450.
\bibitem{wa} M.E. Watkins, On the action of non-abelian groups on graphs, J. Comb. Theory 11 (1971) 95-104.
\bibitem{x2} B. Xia, On cubic graphical regular representations of finite simple groups, J. Combin. Theory Ser. B 141 (2020), 1-30.
\bibitem{xu1}S. Xu, D. Wong, C. Zhang, et al, On oriented $ m $-semiregular representations of finite groups about valency three, arxiv preprint arxiv:2507.15405, 2025.
\bibitem{xu2}S. Xu, D. Wong, C. Zhang, et al, The $ m $-partite digraphical representations of valency 3 of finite groups generated by two elements, arxiv preprint arxiv:2503.22980, 2025.

\end{thebibliography}
\end{document}